\documentclass[11pt,reqno]{amsart}
\usepackage[T1]{fontenc}
\usepackage[utf8]{inputenc}
\usepackage{lmodern,amsmath,amssymb,amsthm,mathrsfs,mathtools,booktabs,microtype}
\usepackage[margin=1in]{geometry}
\usepackage{xcolor,cite,needspace}
\usepackage{tikz,float}
\usetikzlibrary{arrows.meta}
\usepackage[colorlinks=true,linkcolor=blue!45!black,citecolor=blue!50!black,urlcolor=blue!60!black]{hyperref}
\hypersetup{pdftitle={Improved bounds for the binary Chowla conjecture},pdfauthor={Jared Duker Lichtman}}
\allowdisplaybreaks[2]
\numberwithin{equation}{section}
\theoremstyle{plain}
\newtheorem{theorem}{Theorem}[section]
\newtheorem{proposition}[theorem]{Proposition}
\newtheorem{lemma}[theorem]{Lemma}
\newtheorem{corollary}[theorem]{Corollary}
\theoremstyle{definition}

\theoremstyle{remark}
\newtheorem{remark}[theorem]{Remark}

\newcommand{\eps}{\varepsilon}

\newcommand{\proofstep}[1]{\par\smallskip\noindent\textit{#1}\enspace\ignorespaces}

\renewcommand{\leq}{\leqslant}
\renewcommand{\geq}{\geqslant}
\renewcommand{\le}{\leqslant}
\renewcommand{\ge}{\geqslant}
\providecommand{\cref}[1]{\eqref{#1}}
\usepackage{iftex}
\ifPDFTeX
  \DeclareUnicodeCharacter{00A0}{\ifmmode\else\space\fi}
\else
  \catcode"00A0=10\relax
\fi

\begin{document}
\title[Improved bounds for the binary Chowla conjecture]
{Improved bounds for the binary Chowla conjecture}
\author{Jared Duker Lichtman}
\begin{abstract}
For the M\"obius function $\mu$ and every fixed $0 < \eps < 1$, we prove
\begin{align*}
  \sum_{n \le x} \mu(n) \mu(n + 1) \ll_{\eps} \frac{x}{(\log x)^{1 - \eps}}.
\end{align*}
This improves on OpenAI's recent bound of $x/(\log x)^c$, where $c > 0$ is a tiny absolute constant.

The proof echoes the graph expansion framework of Helfgott, Radziwi\l\l {}, and Pilatte on the logarithmic Chowla conjecture, earlier Fourier analytic methods of Matom\"aki, Radziwi\l\l {}, and Tao on the Chowla conjecture on average, as well as work of Ter\"av\"ainen and the author for M\"obius on shifted primes. This also yields an improved result,
\begin{align*}
  \sum_{p \le x} \mu(p + 2) \ll_{\eps} \frac{\pi(x)}{(\log x)^{1/2 - \eps}}.
\end{align*}
Specifically, we refine the determinant approach from the Poisson--Dirichlet law for prime predecessors. We replace the averaged prime weights of L\"u--Wang--Wu with an exact partition of prime intervals. Once sieving out exceptional integers, this enables uniform cancellation with arbitrary log power savings.
\end{abstract}
\maketitle

\section{Introduction}  \label{sec:introduction}
The classical M\"obius function $\mu$ given by $\mu(n) = (-1)^r$ for squarefree $n$ with $r$ prime factors, and $\mu(n) = 0$ otherwise. The prime number theorem is equivalent to $\sum_{n \le X}\mu(n) = o(X)$, while the Riemann hypothesis is equivalent to
\begin{align*}
  \sum_{n \le X}\mu(n) \ll_{\eps} X^{1/2 + \eps}.
\end{align*}
For correlations with $\mu$, the well-known Chowla conjecture~\cite{Chowla} asserts that, for any $k \geq 1$ and distinct shifts $h_1, \ldots, h_{k}$, one has
\begin{equation} \label{conj:chowla}
  \sum_{n \leq X} \mu(n + h_1) \mu(n + h_2) \cdots \mu(n + h_k) = o(X).
\end{equation}
This can be regarded as a multiplicative analogue of the Hardy--Littlewood prime $k$-tuple conjecture, which predicts an asymptotic formula for such correlations. In fact, recently a Hardy--Littlewood--Chowla conjecture was formulated for hybrid correlations~\cite{Lichtman2022,LichtmanTeravainen2022}, including the shifted prime problem for $\sum_{p \le X} \mu(p + 2)$. Note that the study of $\mu$ is essentially interchangeable with that of the Liouville function $\lambda$.

Chowla's conjecture was recently resolved in the case $k = 2$ by OpenAI \cite{OpenAIChowla2026}, and remains open for all $k \geq 3$. Prior to this, remarkable progress had been made over the preceding decade on weaker variants of Chowla's conjecture.

In 2015, Matom\"aki, Radziwi\l\l {} and Tao proved that \cref{conj:chowla} holds \emph{on average} over $h_1, \ldots, h_k$, for every fixed $k \geq 2$\cite{MRT}. A key ingredient in their proof was the landmark work of Matom\"aki and Radziwi\l\l~\cite{MR} on averages of multiplicative functions in short intervals. Soon after, Tao proved the weaker logarithmic version of Chowla's conjecture for $k = 2$\cite{Tao2016}, showing
\begin{equation} \label{eq:logchowla}
  \sum_{n \leq X} \frac{1}{n} \mu(n) \mu(n + 1) = o(\log X).
\end{equation}
The logarithmic version of \cref{conj:chowla} was later proved for all odd $k \geq 3$ by Tao and Ter\"av\"ainen~\cite{TTduke,taoODD}. For even $k \geq 4$, the logarithmic Chowla conjecture remains open. 

In 2021, Helfgott and Radziwi\l\l~\cite{HR2021} proved a remarkable quantitative refinement,
\begin{equation}
  \label{eq:HRchowla} \sum_{n \leq X} \frac{1}{n} \mu(n) \mu(n + 1) \ll \frac{\log X}{(\log \log X)^{1/2}}.
\end{equation}
Their combinatorial proof was quite different, studying the expansion of a certain weighted graph defined in terms of divisibility by small primes. Also see the exposition by Helfgott~\cite{explanation}. By enlarging the graph steps from primes to products of primes, Pilatte \cite{Pilatte} managed to prove the first logarithmic savings, showing for an absolute $c > 0$,
\begin{align*}
\sum_{n \leq X} \frac{1}{n} \mu(n) \mu(n + 1) \ll (\log X)^{1 - c}.
\end{align*}

Recently, OpenAI \cite{OpenAIChowla2026} found a clever argument to convert the logarithmic version to the original Chowla problem \eqref{conj:chowla} when $k = 2$, obtaining
\begin{align*}
  \sum_{n \le X} \mu(n) \mu(n + 1) \ll \frac{X}{(\log X)^c},
\end{align*}
where $c > 0$ is Pilatte's same exponent \cite{Pilatte}. A limit of optimizing their method suggests $c \approx 10^{-5}$.

The main result of this article is to obtain an exponent $c = 1 - \eps$, for any $\eps > 0$.

\begin{theorem}  \label{thm:chowla}
Let $0 < \eps < 1$. Then for all $X > 1$, we have
\begin{equation}  \label{eq:main-chowla}
  \sum_{n \le X} \mu(n) \mu(n + 1) \ll_{\eps} \frac{X}{(\log X)^{1 - \eps}}.
\end{equation}
\end{theorem}

The proof techniques apply equally with the Liouville function $\lambda$ in place of $\mu$, and with an arbitrary fixed shift $n + h$ with $h \neq0$. We proceed with the case $h = 1$ for notational convenience, as in \cite{Pilatte}. Moreover, the same proof method also establishes a companion result for M\"obius on shifted primes, with savings $(\log X)^{1/2 - \eps}$.

\begin{theorem}  \label{thm:shifted}
Let $0 < \eps < 1$. Then for all $X > 1$, we have
\begin{equation}  \label{eq:shifted}
  \sum_{p \le X} \mu(p + 2) \ll_{\eps} \frac{\pi(X)}{(\log X)^{1/2 - \eps}}.
\end{equation}
\end{theorem}

This gives a quantitative refinement of L\"u, Wang, and Wu \cite{LWW} with the qualitative bound $\sum_{p \le X} \mu(p + 2) = o(\pi(X))$. We also note prior work of the author and Ter\"av\"ainen on averaged forms of the shifted prime problem~\cite{Lichtman2022,LichtmanTeravainen2022}. Specifically, work from the author's thesis (related to, e.g., \cite[Lemma 6.1]{Lichtman2022}) implied that results bounding correlations with shifted primes could be adapted and, in fact, amplified into Chowla. Thus, the savings in Theorem \ref{thm:shifted} may be transferred and effectively doubled (for the two factors of $\mu$) into that of Theorem \ref{thm:chowla}.

Another consequence of refining OpenAI's estimates, we quantify the remainder term for smooth shifted primes~\cite[Theorem~3.1 and (1.3)]{PD}, originally conjectured by Pomerance~\cite[p.~86]{Pomerance1980} (also see Granville~\cite[\S5.3, (5.1)]{GranvilleSmooth}).

\begin{theorem}\label{thm:smooth-prime-predecessors}
Let $\rho$ denote the Dickman--de Bruijn function. For fixed $u\ge 1$, we have
\[
  \#\{p\le x : P^+(p-1)\le x^{1/u}\}
  = \rho(u)\pi(x)
    + O_{u,\varepsilon}\!\left(\frac{\pi(x)}{(\log x)^{1/2-\varepsilon}}\right).
\]
\end{theorem}

The same quantitative rate of convergence applies to \cite[Theorem~1.1]{PD} for the joint Poisson--Dirichlet law for prime predecessors, conjectured by Ford, Konyagin, and Luca~\cite[Conjecture~5]{FKL}.

\subsection*{Acknowledgements and AI statement}

The author would like to thank Harald Helfgott for encouragement, helpful comments and feedback. The author used ChatGPT-6 Astra for assistance writing up the article, checking, and transcribing routine adaptations of arguments. The author takes full responsibility for the correctness of the proof.

This manuscript began as an initial effort to digest and distill the multiple recent 70+ page papers by OpenAI and L\"u--Wang--Wu. The author was also encouraged by a recent clarifying exposition of Lamzouri~\cite{Lamzouri2026} on Anthropic's proof of simplicity of zeta zeros~\cite{AlpogeFurman2026}.

Their papers obtained striking results, but appeared to the author quite difficult to read. The author found it much more helpful to ask questions and discuss extensively with ChatGPT-6, rather than read the primary materials themselves.

From increasingly pointed questioning of ChatGPT-6, the author was surprised to realize that the precise cancellation mechanism possessed extra uniformity that was not being exploited in \cite{LWW}, made available since only required one of the variables in the Type II estimate to be rough (c.f. Proposition \ref{rev:uniform-typeii} below).

The author's principal contribution lies in the following two items: 
\begin{itemize}
\item[(i)] To devise the more precise typical/exceptional integer decomposition via prime intervals, rather than the averaged overlapping weights of L\"u--Wang--Wu \cite{LWW}. In fact, the author views this decomposition closer in spirit to the intervals arising in Matomaki--Radziwi\l\l's breakthrough work \cite{MR}. This ultimately facilitated the improvement in Theorem \ref{thm:shifted} over \cite[Theorem 1]{LWW}
    
\item[(ii)] To transfer the typical integer estimate from shifted primes to Theorem \ref{thm:chowla}. Specifically, work from the author's thesis (related to, e.g., \cite[Lemma 6.1]{Lichtman2022}) implied that results bounding correlations with shifted primes could be adapted and, in fact, amplified into Chowla. Thus, the savings in Theorem \ref{thm:shifted} may be transferred and effectively doubled (for the two factors of $\lambda$) into that of Theorem \ref{thm:chowla}.
\end{itemize}

Despite the considerable computational resources devoted to the problem by OpenAI and others, the author found it heartening that human expertise may still discover further `meat left on the bone', to extract for improvement in Theorem \ref{thm:chowla}.

\section{Outline of the proof}  \label{sec:outline}

Choose $0 < \alpha < \beta < 1/2$ and $k\ge1$ large. Writing $\mathcal{L} = \log X$, we split the interval $[\exp(\mathcal{L}^\alpha),\, \exp(\mathcal{L}^\beta)]$ into consecutive intervals of the form $I = (y,y^2]$. Note there are $B$ such intervals $I$, where $2^B = \mathcal{L}^{\beta - \alpha}$, i.e.
\begin{align*}
  B = \frac{\beta - \alpha}{\log2} \log \mathcal{L} + O(1).
\end{align*}
Also, the density of integers without a prime factor in $I$ is $\prod_{p \in I} (1 - 1/p) \sim 1/2$.

Call an interval successful for $n$ if it contains a prime divisor of $n$. Call $n$ \emph{typical} if at least $3k$ intervals are successful, and \emph{exceptional} otherwise. Let $T,E$ be their indicators, so $1=T(n)+E(n)$. From the identity,
\begin{align*}
  1 = T(n) + E(n)T(n + 1) + E(n)E(n + 1),
\end{align*}
we multiply by $\mu(n) \mu(n + 1)$ and sum over $n \sim X$. Thus by the triangle inequality,
\begin{align}
  \bigg| \sum_{n \sim X} \mu(n) \mu(n + 1) \bigg| \le \bigg| \sum_{n \sim X} \mu(n)T(n) \mu(n + 1) \bigg| & + \bigg| \sum_{n \sim X} \mu(n)E(n) \mu(n + 1)T(n + 1) \bigg| \label{uniform:intro-two-residual} \\
  & + \sum_{n \sim x}E(n)E(n + 1). \label{uniform:intro-two-residual2}
\end{align}

Recalling the density $\prod_{p \in I} (1 - 1/p) \sim 1/2$, we will show that the exceptional density is at most $2^{-B}B^{O_k(1)} = \mathcal{L}^{-(\beta - \alpha) + o(1)}$. We make this sieve argument precise in Section~\ref{sec:dense-reconstruction} to give the bound
\begin{align*}
  \sum_{n \sim X}E(n) \le X \mathcal{L}^{-(\beta - \alpha) + o(1)}.
\end{align*}
For consecutive exceptional integers, an analogous sieve argument in Lemma~\ref{uniform:joint-residual} will give
\begin{equation}  \label{eq:intro-residual}
  \sum_{n \sim X}E(n)E(n + 1) \le x \mathcal{L}^{-2(\beta - \alpha) + o(1)}.
\end{equation}
This handles the sum in \eqref{uniform:intro-two-residual2}. We shall prove stronger estimates when $n$ or $n + 1$ is typical,
\begin{align} \label{eq:typicals}
  \bigg| \sum_{n \sim X} \mu(n)T(n) \mu(n + 1) \bigg|, \; \bigg| \sum_{n \sim X} \mu(n)E(n) \mu(n + 1)T(n + 1) \bigg| \ll_A X\mathcal{L}^{-A}.
\end{align}
Substituting \eqref{eq:intro-residual}--\eqref{eq:typicals} back into \eqref{uniform:intro-two-residual} gives
\begin{equation}  \label{eq:dyadic-goal}
  \bigg| \sum_{n \sim X} \mu(n) \mu(n + 1) \bigg| \le X \mathcal{L}^{-2(\beta - \alpha) + o(1)} + O_A(X \mathcal{L}^{-A}).
\end{equation}
Taking $\alpha = \eps/8$ and $\beta = 1/2 - \eps/8$ gives $2(\beta - \alpha) = 1 - \eps/2$. Summing \eqref{eq:dyadic-goal} dyadically then gives Theorem~\ref{thm:chowla}.

\subsection{The typical-integer estimate}
For the estimate \eqref{eq:typicals}, choose $\beta < \zeta < 1$ and set $z = \exp(\mathcal{L}^ \zeta)$. Write $s_n$ for the $z$-smooth part of $n$. Since all available primes lie below $z$, both $T(n)$ and $E(n)$ depend only on $s_n$. That is, $T(n) = T(s_n)$ and $E(n) = E(s_n)$.

Uniformly for any 1-bounded weight $\mathcal{U}$, in Corollary~\ref{uniform:marked-pair} we shall prove a strong \emph{typical-integer estimate},
\begin{equation}  \label{eq:intro-completed}
  \sum_{n \sim X} \mu(n) \mathcal{U}(s_n) \mu(n + 1)T(n + 1) \ll_A X \mathcal{L}^{-A},
\end{equation}
The same estimate holds if we interchange the roles of $n$ and $n + 1$. Hence \eqref{eq:intro-completed} establishes the desired estimate \eqref{eq:typicals}, for the choices of $\mathcal{U}(s_n) = 1$ and $\mathcal{U}(s_n) = E(s_n) = E(n)$, respectively. (For the former choice, we interchange the roles of $n$ and $n + 1$.) Indeed, the key feature of the estimate \eqref{eq:intro-completed} is that $\mathcal{U}$ only needs the $z$-smooth part $s_n$ as input.

Partition typical integers by their first $3k$ successful intervals, with increasing indices $J = \{j_i:0 \le i < 3k\}$. Their pattern indicator $\phi_J$ gives $T = \sum_J \phi_J$. Take $J^+ = \{j_{3i}:0 \le i < k\}$ as active and set $J^- = J \setminus J^+$. Write $\phi_J = \phi^+\phi^-$, where $\phi^+$ records active hits and $\phi^-$ the passive hits and forced misses. Successive active logarithmic scales have ratio at least $8$. There are $O_k(B^{3k}) = \mathcal{L}^{o(1)}$ patterns, so it suffices to prove \eqref{eq:intro-completed} for each $\phi_J$. Section~\ref{sec:arithmetic-reduction} decomposes $\mu(n)\mathcal{U}(s_n)$ into Type~I and Type~II sums.

At this stage, the Vaughan decomposition itself is standard, yet applying it to $\mu$ is perhaps unexpected in the context of the Chowla problem. The Type~I estimate follows from a standard Bombieri--Vinogradov estimate for $\mu$. Both Type I and II sums end up bounded by $X\mathcal{L}^{-A}$ with arbitrary log savings.

\Needspace{15\baselineskip}

\begin{figure}[H]
\centering
\hypersetup{linkcolor=blue}
\newcommand{\proofref}[2]{#1~\ref{#2}}
\newcommand{\proofdesc}[1]{%
  \parbox{3.00cm}{\centering\fontsize{8.5}{10}\selectfont #1}}
\resizebox{\linewidth}{!}{%
\begin{tikzpicture}[
  result/.style={
    draw=black!70,
    line width=0.45pt,
    rounded corners=3pt,
    fill=white,
    text width=3.00cm,
    align=center,
    inner xsep=4pt,
    inner ysep=4pt,
    minimum height=0.95cm,
    font=\fontsize{10}{12}\selectfont
  },
  conclusion/.style={result,double,double distance=1.1pt},
  implication/.style={
    -{Latex[length=1.7mm,width=1.1mm]},
    draw=black!85,
    line width=0.45pt,
    shorten <=1.2pt,
    shorten >=1.2pt
  }
]

\node[result] (bv) at (4,0)
  {\proofref{Lemma}{lem:mobius-BV}\\[2pt]
   \proofdesc{Weighted\\Bombieri--Vinogradov}};
\node[result] (mobius) at (8,0)
  {\proofref{Proposition}{prop:mel}\\[2pt]
   \proofdesc{M\"obius cancellation\\Dirichlet polynomials}};
\node[result] (moment) at (12,0)
  {\proofref{Proposition}{prop:par}\\[2pt]
   \proofdesc{Determinant\\graph expansion}};

\node[result] (typei) at (4,-1.85)
  {\proofref{Proposition}{rev:mobius-typeI}\\[2pt]
   \proofdesc{Type I sums}};
\node[result] (major) at (8,-1.85)
  {\proofref{Proposition}{prop:maj}\\[2pt]
   \proofdesc{Major arcs}};
\node[result] (minor) at (12,-1.85)
  {\proofref{Proposition}{prop:min}\\[2pt]
   \proofdesc{Minor arcs}};

\node[result] (typeii) at (8,-3.60)
  {\proofref{Proposition}{rev:uniform-typeii}\\[2pt]
   \proofdesc{Type II estimate}};

\node[result] (arithmetic) at (4,-3.60)
  {\proofref{Proposition}{prop:vau}\\[2pt]
   \proofdesc{Vaughan decomposition\\Type I and II}};

\node[result] (ramare) at (12,-3.60)
  {\proofref{Lemma}{lem:ramare-identity}\\[2pt]
   \proofdesc{Ramar\'e identity}};

\node[result] (partition) at (0,-3.60)
  {\proofref{Identity}{gr:partition}\\[2pt]
   \proofdesc{Typical / exceptional\\decomposition}};

\node[result] (typical) at (7.5,-5.35)
  {\proofref{Corollary}{uniform:marked-pair}\\[2pt]
   \proofdesc{Typical integer estimate}};
\node[result] (sieve) at (0,-5.35)
  {\proofref{Lemma}{uniform:joint-residual}\\[2pt]
   \proofdesc{Exceptional integer\\sieve bound}};

\node[conclusion] (theorem) at (4,-7.10)
  {\proofref{Theorem}{thm:chowla}\\[2pt]
   \proofdesc{Binary Chowla bound}};

\draw[implication] (bv) -- (typei);
\draw[implication] (mobius) -- (major);
\draw[implication] (moment.south) -- (minor.north);

\draw[implication] (major.south) -- (typeii.north);
\draw[implication] (minor.south) -- (typeii.north east);
\draw[implication] (ramare.west) -- (typeii.east);

\draw[implication] (typei.south east)
  to[out=-72,in=90] (6,-3.15)
  -- (6,-4.22)
  .. controls (6,-4.55) and (6.80,-4.44)
  .. ([xshift=-7mm]typical.north);
\draw[implication] (arithmetic.south east) -- (typical.north west);
\draw[implication] (typeii.south)
  -- ([xshift=3mm]typical.north);
\draw[implication] (partition.south east)
  to[out=-25,in=180,looseness=0.8] (typical.west);
\draw[implication] (partition) -- (sieve);

\draw[implication] (typical.south) -- (theorem.north east);
\draw[implication] (sieve.south) -- (theorem.north west);

\end{tikzpicture}
}
\caption{The dependency structure in the proof of \proofref{Theorem}{thm:chowla}.}
\label{fig:proof-outline}
\end{figure}
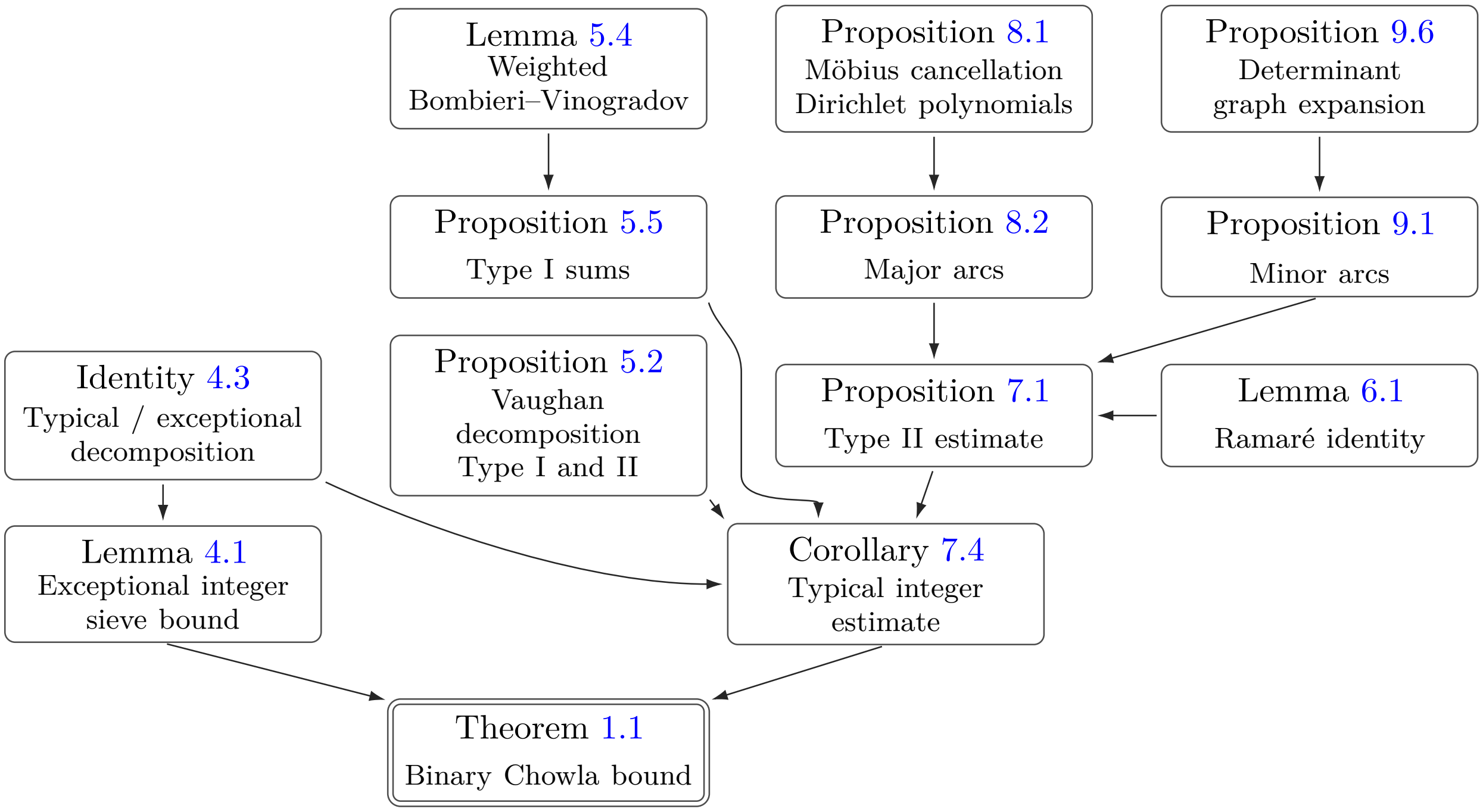

The main result considers Type~II sums of the form
\begin{equation}  \label{eq:intro-typeii}
  \sum_{\substack{m \asymp M \\P^-(m) > z}} \sum_{\substack{n \asymp N \\mn \in I}} \mu(m)\mathcal{B}(n) \mu(mn + 1) \phi(mn + 1).
\end{equation}
Here $\mathcal{B}(n)$ is a certain divisor-bounded coefficient, which contains $\mu(s_n) \mathcal{U}(s_n)$. A key feature of the sum \eqref{eq:intro-typeii} is that $m$ is assumed to be $z$-rough, while $n$ is unrestricted.

We extract $k$ primes from intervals in the pattern $\phi$, via the Ramar\'e identity. Grouping cofactors and applying Cauchy--Schwarz leads to bilinear sums indexed by a determinant equation,
\begin{equation}  \label{eq:intro-determinant}
 \det\begin{pmatrix} bm & ar \\ s & n \end{pmatrix} =bmn-ars
  = a - b.
\end{equation}
Here $a,b$ are products of the $k$ chosen primes, and $r,s$ are copies of the variables $m,n$, respectively. The bilinear sum itself look roughly like
\begin{align*}
Q=\sum_{\substack{m,n,r,s\\ P^-(mr) > z}}\sum_{\substack{bmn - ars = a - b\\ (a,b) = 1}} \alpha_m \overline{\alpha_r} \beta_n \overline{\beta_s}
\end{align*}
for certain divisor-bounded coefficients $\alpha_m,\beta_n$. In Section~\ref{sec:comparison}, we shall decouple the variables from the determinant equation \eqref{eq:intro-determinant}, via the Fourier identity with $t=bmn - ars$,
\begin{align*}
\mathbf{1}_{t=a - b} = \int_0^1 e\big(\theta(t - a + b)\big)\,d\theta.
\end{align*}
Then we split the resulting bilinear form into major arcs $\frak{M}$ and minor arcs $\frak{m}$.

To prove the major arc bound in Proposition~\ref{prop:maj}, we use a $z$-rough adaptation of a Dirichlet polynomial estimate of Matom\"aki--Shao--Tao--Ter\"av\"ainen \cite{MSTT}. The cancellation ultimately derives from the Vinogradov--Korobov zero-free region for Dirichlet $L$-functions.

\subsection{Minor arcs}

To prove the minor arc bound in Proposition~\ref{prop:min}, we combine Fourier orthogonality with a count of live primes becoming active again. This refines the determinant framework of OpenAI \cite{PD} and L\"u--Wang--Wu \cite{LWW} by only requiring rough support for the first coefficient sequence, and allows $\mathcal{B}(n)$ to keep the smooth-part weight. Notably, neither source of saving requires cancellation in the coefficients' signs.

The minor arcs proof is then a technical adaption of L\"u--Wang--Wu \cite{LWW}, which we outline for completeness. The overarching idea is to average over prime-labelled paths, while tracking divisibility conditions. Indeed, this echoes Helfgott--Radziwi{\l}{\l} \cite[Sections~1.5 and~5]{HR2021} and Pilatte \cite[Proposition~5.3]{Pilatte}.

We shall bound the sum over minor arcs as $Q^{\rm min}  \ll_A X\mathcal{L}^{-A}$, which by suppressing technicalities, looks like
\begin{align*}
  Q^{\rm min}\simeq\sum_{\substack{m,n,r,s\\P^-(mr)>z}}\sum_{a,b\in\mathcal A}\alpha_m\overline{\alpha_r}\beta_n\overline{\beta_s}\int_{\mathfrak m}e\bigl(\theta(bmn-ars-a+b)\bigr)\,d\theta.
\end{align*}

We express the minor arc sum $Q^{\rm min}$ inner products $\langle f,\mathsf T f\rangle$ for a certain matrix of phases $\mathsf T$ with a coefficient vector $f$, which we now explain.

Following \cite{LWW,PD}, adjoin $J$ auxiliary primes from each interval $\mathcal P_i$, and let $d$ be their product. Set
\begin{align*}
  \mathsf{v}=\binom{bd m}{s},\qquad
  \mathsf{v}'=\binom{ad r}{n},\qquad
  \frac{\det(\mathsf{v},\mathsf{v}')}{d}=bmn-ars.
\end{align*}

Fix the scale $S\asymp bd, ad$, and let $B^*$ consist of the primitive vectors $\mathsf{v}=\binom{\mathsf x}{\mathsf y}$ with $\mathsf x\asymp SM$ and $\mathsf y\asymp N$. A state $(\mathsf{v},c)$ consists of $\mathsf{v}\in B^*$ and a squarefree divisor $c\mid\mathsf x$ formed from $J+1$ primes in each $\mathcal P_i$. The transition matrix $\mathsf T$ averages uniformly over the $(J+1)^k$ choices of $b\mid c$ formed by choosing one prime from each interval. It replaces $b$ by $a\in\mathcal A$ coprime to $c$, keeping $d=c/b$ fixed. Suppressing cutoffs and normalizations, $\mathsf T$ is essentially given by
\begin{align*}
  \mathsf T_{(\mathsf{v},bd),(\mathsf{v}',ad)}
  \simeq \int_{\mathfrak m}
  e\!\left(\theta\left(\frac{\det(\mathsf{v},\mathsf{v}')}{d}+b-a\right)\right)\,d\theta,
  \quad \textrm{ supported on vectors with} \ \ |\det(\mathsf{v},\mathsf{v}')|\ll S.
\end{align*}
The cutoff $|\det(\mathsf{v},\mathsf{v}')|\ll S$ effectively makes $\mathsf T$ into a $O(S)$-near-diagonal matrix. This is a 2-dimensional analog of the near-diagonal matrix in Helfgott--Radziwi{\l}{\l} and Pilatte, which is a central mechanism for spectral expansion.

For each state, define $f(\mathsf{v},c)=\overline{\alpha_{\mathsf x/c}}\,\beta_{\mathsf y}$ so that $\overline{f(\mathsf{v})}f(\mathsf{v}')=\alpha_m\overline{\alpha_r}\beta_n\overline{\beta_s}$. This essentially gives
\begin{align*}
Q_S^{\rm min}
  \simeq \langle f,\mathsf T f\rangle.
\end{align*}

We project $\mathsf T$ onto `good' states at every path step, which induces the self-adjoint operator $\mathsf A=\mathsf G\mathsf T\mathsf G$. Here `good' $\mathsf G$ is just a mild condition on states for the lattice box's $B^*$ length and for well-spaced choices of new primes in $a$, which become important for extracting savings. Then, at scale $S$, we essentially obtain
\begin{align*}
  Q_S^{\rm min}
  \simeq \langle f,\mathsf T f\rangle
  \simeq \langle f,\mathsf A f\rangle.
\end{align*}
Lemmas~\ref{lem:transition-reconstruction} and~\ref{lem:operator-reconstruction} give the precise formulas.

Then by Cauchy--Schwarz and H\"older's inequality for the positive operator $\mathsf A\mathsf A^*=\mathsf A^2$,
\begin{align*}
  |\langle f,\mathsf A f\rangle| \le \|f\|\,\|\mathsf A^*f\| 
  &\le \|f\|\cdot \|f\|^{1-1/R}\left\langle f,(\mathsf A\mathsf A^*)^R f\right\rangle^{1/(2R)}\\
  & = \|f\|^{2-1/R}\cdot \|A^R f\|^{1/R}.
\end{align*}
Here the coefficients have norms roughly $ \|f\| \ll \sqrt{XS}$. The main task is to bound the high moment,
\begin{align}\label{eq:highmom}
\|A^R f\|^2 \ll XRS^2  \,\mathcal L^{-2RE}.
\end{align}
If so, then plugging back in gives the desired bound,
\begin{align}\label{eq:Qminfin}
Q_S^{\rm min} \simeq |\langle f,\mathsf A f\rangle|
&\ll (\sqrt{XS})^{2-1/R} (XRS^2\,\mathcal L^{-2RE})^{1/(2R)} \nonumber\\
&= XS (RS)^{1/(2R)}\,\mathcal L^{-E} \ll XS \,\mathcal L^{-E},
\end{align}
provided $S^{1/(2R)} \ll 1$ (note $R^{1/(2R)} \ll 1$). That is, \eqref{eq:Qminfin} follows from \eqref{eq:highmom} and $\log S \ll R$.

To prove \eqref{eq:highmom}, consider the orthogonal basis of indicators $u_{\mathsf{v}}(\mathsf{v}',c)=\mathbf1_{\mathsf{v}'=\mathsf{v}}$. Note $u_{\mathsf{v}}$ are unnormalized. Since $f(\mathsf{v},c)=f(\mathsf{v})$, we have $f=\sum_{\mathsf{v}\in B^*}f(\mathsf{v})u_{\mathsf{v}}$. For simplicity, pretend $|f|\ll 1$ hereafter. 
Then by a short argument in Lemma~\ref{lem:block-interpolation}, $\langle u_{\mathsf{v}},\mathsf A^{2R}u_{\mathsf{v}'}\rangle=0$ unless $|\det(\mathsf{v},\mathsf{v}')|\ll RS$, and that
\begin{align}\label{eq:detboxbd}
  \#\big\{\mathsf{v}'\in B^*:|\det(\mathsf{v},\mathsf{v}')|\ll RS\big\}\ll RS.
\end{align}
Thus Cauchy--Schwarz and \eqref{eq:detboxbd} give
\begin{align}
  \|A^R f\|^2
  &=\sum_{\mathsf{v},\mathsf{v}'\in B^*}\overline{f(\mathsf{v})}f(\mathsf{v}')\left\langle\mathsf A^Ru_{\mathsf{v}},\mathsf A^Ru_{\mathsf{v}'}\right\rangle \nonumber\\
  &\ll\sum_{\substack{\mathsf{v},\mathsf{v}'\in B^*\\|\det(\mathsf{v},\mathsf{v}')|\ll RS}}\left(\|\mathsf A^Ru_{\mathsf{v}}\|^2+\|\mathsf A^Ru_{\mathsf{v}'}\|^2\right) \
  \ll \ RS\,\mathcal T_R, \label{eq:momtoTR}
\end{align}
for the (weighted) trace $\mathcal T_R=\sum_{\mathsf{v}\in B^*}\|\mathsf A^Ru_{\mathsf{v}}\|^2$. Proposition~\ref{prop:par} shall give
\begin{align}\label{eq:TRbound}
  \mathcal T_R \ll_{E,k} XS\,\mathcal L^{-2RE}.
\end{align}
Combining \eqref{eq:momtoTR} and \eqref{eq:TRbound} gives
\eqref{eq:highmom}, hence the bound $ Q_S^{\rm min}\ll_{E,k} XS\mathcal L^{-E}$ in \eqref{eq:Qminfin}. This leads to the final minor arcs $|Q^{\rm min}|\ll_{E,k}X\mathcal L^{-E}$ as desired. 

Now to justify the bound for \(\mathcal T_R\) in \eqref{eq:TRbound}, expand \(\mathsf A^{2R}\) as a sum over length-$2R$ paths. Since \(u_{\mathsf{v}}\) fixes only the vector, the vectors trace a loop based at \(\mathsf{v}\). So collecting transition data, normalizations, and cutoffs into weights $W$,
\[
\begin{aligned}
\mathcal T_R
=\sum_{\mathsf{v}\in B^*}\langle u_{\mathsf{v}},\mathsf A^{2R}u_{\mathsf{v}}\rangle
=\sum_{\substack{\mathsf{v}_0,\ldots,\mathsf{v}_{2R}\in B^*\\\mathsf{v}_{2R}=\mathsf{v}_0}}
  \sum_{\substack{c_0,\ldots,c_{2R}\\
                  c_j\mid x_j\ \text{for all }j}}
  W\!\left((\mathsf{v}_j,c_j)_{j=0}^{2R}\right),
\end{aligned}
\]
where \(x_j\) is the first coordinate of \(\mathsf{v}_j\).

Fix the base point \(\mathsf{v}_0=\binom{x_0}{y_0}\). The coordinates $x_0$, $y_0$ are coprime so there exists a unimodular matrix \(G\) such that $\mathsf{v}_0=G\binom{1}{0}$. Writing \(\mathsf{v}_j=G\mathsf{w}_j\), we obtain a transformed loop $(\mathsf{w}_j)_j$ with standardized base point \(\mathsf{w}_0=\mathsf{w}_{2R}=\binom{1}{0}\). Since \(\mathsf{v}_j=\binom{x_j}{y_j}\), the condition \(p\mid x_j\) says that \(\mathsf{v}_j\bmod p\) always lies on the vertical line \(\mathbb F_p\binom{0}{1}\). Then since $G$ is unimodular, \(G\) is invertible and denote the lines $\ell_0(\mathsf{v_0})=G^{-1}\!\left(\mathbb Z\binom{0}{1}\right)$ and \(\ell_p(\mathsf{v_0})=\bar G^{-1}\!\left(\mathbb F_p\binom{0}{1}\right)\). Hence for all $j$,
\[
\begin{aligned}
p\mid x_j
\ \Longleftrightarrow \ \mathsf{v}_j\bmod p\in\mathbb F_p\tbinom{0}{1}
\ \Longleftrightarrow \ \bar G\bar{\mathsf{w}}_j\in\mathbb F_p\tbinom{0}{1}
\ \Longleftrightarrow \ \bar{\mathsf{w}}_j\in\bar G^{-1}\big(\mathbb F_p\tbinom{0}{1}\big)=\ell_p(\mathsf{v_0}).
\end{aligned}
\]
That is, OpenAI has constructed a rather striking mechanism: given a base point $\mathsf{v}_0$ there is a global (geometric) line $\ell_0=\ell_0(\mathsf{v_0})$ that encodes the local arithmetic for the path, i.e. $p\mid x_j$ if and only if $\mathsf{w}_j$ lies on the reduction $\ell_p$ of this line.

Now fix the relative path \(\gamma=(\mathsf{w}_j,c_j)_{j=0}^{2R}\). Since \(\mathsf{v}_0\) determines \(G\), it determines every \(\mathsf{v}_j=G\mathsf{w}_j\). Thus, for fixed \(\gamma\), the weight \(W=W_\gamma(\mathsf{v}_0)\) depends only on \(\mathsf{v}_0\). Put
\[
J_p=\{j:p\mid c_j\},
\qquad
\mathcal P_\gamma=\{p:J_p\ne\varnothing\}.
\]
At the first index \(j\in J_p\), the vector \(\bar{\mathsf{w}}_j\) is nonzero by primitivity, so it spans a line \(L_p\). The path can contribute only if \(\bar{\mathsf{w}}_j\in L_p\) at every later index in \(J_p\). Write \(\sum^*_{\gamma}\) for the sum over paths satisfying this condition for every \(p\in\mathcal P_\gamma\).

For such a path, the divisibilities \(c_j\mid x_j\) hold precisely when \(\mathsf{v}_0\) gives \(\ell_p(\mathsf{v_0})=L_p\) for every \(p\in\mathcal P_\gamma\). Regrouping by relative path $\gamma$ and then averaging over the base point $P_0$ gives
\[
\mathcal T_R =\sideset{}{^*}\sum_{\gamma}
  \sum_{\substack{\mathsf{v}_0\in B^*\\
        \ell_p(\mathsf{v_0})=L_p\ \forall\,p\in\mathcal P_\gamma}}
  W_\gamma(\mathsf{v}_0) 
  \ \simeq \ \sideset{}{^*}\sum_{\gamma} \bigg(\sum_{\mathsf{v}_0\in B^*}W_\gamma(\mathsf{v}_0) \bigg)
  \prod_{p\in\mathcal P_\gamma}\frac1{p+1}.
\]
Indeed, as the base $\mathsf{v_0}$ varies, $\ell_p(\mathsf{v_0})$ equidistributes among the $p+1$ projective lines in $\mathbb{P}^1(\mathbb{F}_p)$, independently at distinct $p$ by Corollary~\ref{cor:root-slopes}. Thus each prime contributes a factor $1/(p+1)$ at its first appearance. If $p$ returns, the star requires the returning vector to lie on the prescribed line $L_p$.

Consider the transition from $c_t$ to $c_{t+1}$. Write
\begin{align*}
c=c_t=\prod_{i=1}^k\prod_{j=0}^J q_{i,j},\qquad
b=\prod_{i=1}^k q_{i,j_i},\qquad
a=\prod_{i=1}^k p_i.
\end{align*}
Given $c$, there are $(J+1)^k$ choices of the factor $b$. The transition replaces each $q_{i,j_i}$ with an incoming prime $p_i$ from the same interval $\mathcal P_i$. Writing $d=c/b$, we have $c_t=bd$ and $c_{t+1}=ad$.

Link each departure of a prime to its next arrival, if any. A step is linked if it belongs to such a pair, and unlinked otherwise. At step $t$, denote the indicators that every incoming prime is born by $B_t=\prod_{p\mid a}\mathbf 1_{\min J_p=t+1}$, and that every departing prime dies by $D_t=\prod_{p\mid b}\mathbf 1_{\max J_p=t}$. Then we have the identity
\begin{align*}
1&=\underbrace{B_tD_t}_{\text{unlinked}}
  +\underbrace{(1-B_t)}_{\text{a link ends at $t$}}
  +\underbrace{B_t(1-D_t)}_{\text{a link starts (no ends) at $t$}}.
\end{align*}

{\bf Counting on linked steps.}
Reversing the path if necessary, let $p\mid a$ be the largest returning prime. Its earlier appearance prescribes $L_p$. Note that $c$ and $ad/p$ determine $b$ via
\begin{align*}
\gcd\!\left(c,\frac{ad}{p}\right)=d,
\qquad b=\frac{c}{d}.
\end{align*}

For any candidate choice of primes, the line $L_p$ ends up forcing $(ad/p)\rho\bmod 1$ to lie in a fixed window of width $O(1/p)$, where $\rho$ depends only on the initial step state. However, by construction, all prime factors from typical integers are from well-spaced intervals. Hence the products $ad/p$ will also be well-spaced, so only one choice of $ad/p$ can fit the window $\bmod 1$, and hence only one unique $b$ is possible, among the $(J+1)^k$ choices. Averaging over all factors, this gives the final savings
\begin{align*}
(J+1)^{-k} \ll_E \mathcal{L}^{-E},
\end{align*}
setting $J=\mathcal{L}^{\delta}$ with $\delta>0$ small and $k=O(1/\delta)$ large.

{\bf Fourier cancellation on unlinked steps.}
On an unlinked step, the primes in $a$ have never appeared before, and the primes in $b$ will never appear again. As such, no incoming prime imposes a line $L_p$ from the past. Crucially, this allows us to absorb all preceding steps into a weight $W(\mathsf{k},b)$ depending only on the current vector $\mathsf{k}\in\mathbb{Z}^2$ and the departing factor $b$, once the shared data are fixed. Then Cauchy--Schwarz reduces to essentially
\begin{align}\label{eq:indepminor}
  \int_{\mathfrak m}\sum_{\mathsf{l}}
  \left|
  \mathop{\vcenter{\hbox{\scalebox{1.5}{$\mathbb{E}$}}}}\limits_{\mathsf{k}}\, e\big(\theta\det(\mathsf{k},\mathsf{l})\big)
  \sum_{b}W(\mathsf{k},b)e(\theta b)
  \right|^2\,d\theta
  \ll_{E} \,\mathcal L^{-E}\|W\|^2.
\end{align}
Since the phase factors as $e\big(\theta\det(\mathsf{k},\mathsf{l})\big) = e(\theta k_1l_2)\,e(-\theta k_2l_1)$, the integrand separates into bilinear Fourier sums of $e(\alpha mp)$, which take the exact form bounded by Matom\"aki--Radziwi{\l}{\l}--Tao's minor-arc argument in \cite[Section~3]{MRT}. The savings $O(\mathcal{L}^{-E})$ in \eqref{eq:indepminor} ultimately comes from bounding geometric series with Vinogradov's lemma, uniformly over $\theta\in \frak{m}$.

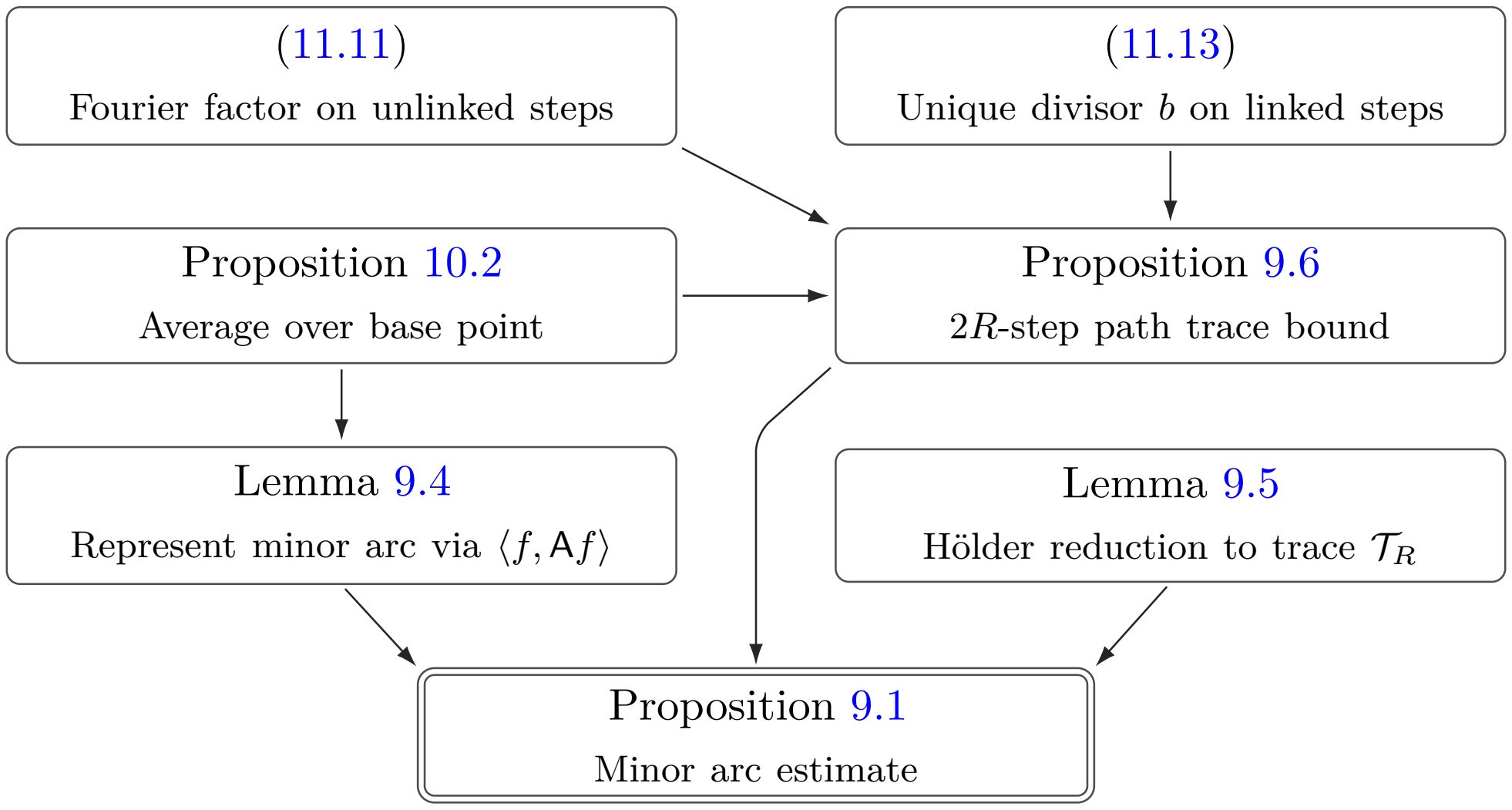
\begin{figure}[H]
\centering
\hypersetup{linkcolor=blue}
\newcommand{\minordesc}[1]{\mbox{\fontsize{8.5}{10}\selectfont #1}}
\resizebox{0.78\linewidth}{!}{%
\begin{tikzpicture}[
  result/.style={
    draw=black!70,
    line width=0.45pt,
    rounded corners=3pt,
    fill=white,
    text width=5.05cm,
    align=center,
    inner xsep=4pt,
    inner ysep=4pt,
    minimum height=0.95cm,
    font=\fontsize{10}{12}\selectfont
  },
  conclusion/.style={result,double,double distance=1.1pt},
  implication/.style={
    -{Latex[length=1.7mm,width=1.1mm]},
    draw=black!85,
    line width=0.45pt,
    shorten <=1.2pt,
    shorten >=1.2pt
  }
]

\node[result] (unlinked) at (-3.30,0)
  {\eqref{eq:cleanbd}\\[2pt]
   \minordesc{Fourier factor on unlinked steps}};

\node[result] (linked) at (3.30,0)
  {\eqref{eq:epsk}\\[2pt]
   \minordesc{Unique divisor $b$ on linked steps}};

\node[result] (roots) at (-3.30,-1.75)
  {Proposition~\ref{prop:root-average}\\[2pt]
   \minordesc{Average over base point}};

\node[result] (moment) at (3.30,-1.75)
  {Proposition~\ref{prop:par}\\[2pt]
   \minordesc{$2R$-step path trace bound}};

\node[result] (operator) at (-3.30,-3.50)
  {Lemma~\ref{lem:operator-reconstruction}\\[2pt]
   \minordesc{Represent minor arc via $\langle f,\mathsf A f\rangle$}};

\node[result] (holder) at (3.30,-3.50)
  {Lemma~\ref{lem:block-interpolation}\\[2pt]
   \minordesc{H\"older reduction to trace $\mathcal T_R$}};

\node[conclusion] (minor) at (0,-5.25)
  {Proposition~\ref{prop:min}\\[2pt]
   \minordesc{Minor arc estimate}};

\draw[implication] (unlinked.south east) -- (moment.north west);
\draw[implication] (linked.south) -- (moment.north);
\draw[implication] (roots.east) -- (moment.west);
\draw[implication] (roots.south) -- (operator.north);
\draw[implication] (operator.south) -- (minor.north west);
\draw[implication,rounded corners=4pt]
  (moment.south west) -- (0,-2.85) -- (minor.north);
\draw[implication] (holder.south) -- (minor.north east);

\end{tikzpicture}%
}
\caption{Proof dependencies for the minor arcs in Proposition~\ref{prop:min}.}
\label{fig:minor-outline}
\end{figure}

Hence both linked and unlinked cases yield a savings factor $O
(\mathcal{L}^{-E})$, at each of the $2R$ steps. Altogether this gives total savings $O(\mathcal{L}^{-2RE})$, so that
\begin{align}
  \mathcal{T}_R\ll_E XS\,\mathcal{L}^{-2RE}.
\end{align}
Together with \eqref{eq:momtoTR}, we obtain \eqref{eq:highmom} to complete the proof of the minor arcs estimate. 

Note that neither minor arc step uses M\"obius cancellation, instead saving via Fourier analysis (for unlinked steps) or counting bounds (for linked steps).

\section{Notation and preliminaries}  \label{sec:prelim}

Fix a sufficiently large $x$ and write $\mathcal{L} = \log x$. Write $n \sim x$ to denote $x < n \le 2x$. When we shift a sum from $n$ to $n + 1$, we keep the prime intervals and cutoffs defined at $x$. Thus a change of variable preserves the original test for typicality. Implied constants may depend on the stated fixed parameters.

We write $*$ for Dirichlet convolution. Let $P^ + (n)$ and $P^-(n)$ be the largest and smallest prime factors, with $P^ + (1) = 1$ and $P^-(1) = \infty$. The function $\tau_j(n)$ counts ordered factorizations of $n$ into $j$ positive integers; $\tau = \tau_2$. We also set $\mathbb T = \mathbb R/ \mathbb Z$, $e(t) = \exp(2 \pi it)$, and write $\|t \|$ for distance to the nearest integer. An interval restriction on a coefficient sequence means that it is set to zero outside that interval.

Fix parameters
\begin{equation}  \label{par:all}
  0 < \alpha < \beta < \tfrac12, \qquad \beta < \zeta < 1, \qquad \jmath > 0, \qquad \beta + \jmath < \varrho < 1 - \beta - \jmath.
\end{equation}
The exponents $\alpha$ and $\beta$ determine what range of prime factors we can inspect. The exponent $\zeta$ determines the roughness cutoff $z$. Later in the minor arc proof, $\jmath$ will set the length of each auxiliary prime list and $\varrho$ will set the path length in the moment expansion. The lower bound on $\varrho$ will allow us to recover the bilinear estimate from the moment. The upper bound on $\varrho$ controls the congruence modulus and the total equidistribution error after summing over paths. We shall ultimately take $\varrho = 1/2$ and $\jmath$ sufficiently small to satisfy both bounds. Set
\begin{equation}  \label{eq:001}
  z = \exp(\mathcal{L}^ \zeta), \qquad R_z(n) = \mathbf1_{P^-(n) > z}, \qquad \nu(n) = \mu(n)R_z(n).
\end{equation}
We call $n$ with $R_z(n) = 1$ \emph{$z$-rough}, and note $R_z(ab) = R_z(a)R_z(b)$. Write $s_n = \prod_{p \le z}p^{v_p(n)}$ for the \emph{$z$-smooth} part of $n$.

An admissible active family consists of a fixed number $k$ of prime intervals
\begin{equation}  \label{rev:admissible-intervals}
  \begin{gathered}
  \mathcal{P}_i = \{p \in \mathbb P:e^{\lambda_i} < p \le e^{2 \lambda_i} \}, \qquad \mathcal{L}^ \alpha \le \lambda_1 < \cdots < \lambda_k \le \mathcal{L}^ \beta, \\
  \lambda_{i + 1} \ge r_{\rm sep} \lambda_i, \qquad r_{\rm sep} > 19/3.
  \end{gathered}
\end{equation}
Define
\begin{align*}
  V_i &= \sum_{p \in \mathcal{P}_i} \frac1p, \qquad V = \prod_{i=1}^k V_i, \qquad a_i = \frac{\log \lambda_i}{\log \mathcal{L}}\\
  p_{\min} &= \min \bigcup_i \mathcal{P}_i, \qquad p_{\max} = \max \bigcup_i \mathcal{P}_i.
\end{align*}
The prime number theorem and partial summation give
\begin{equation}  \label{eq:interval-mass}
  V_i = \log2 + O_A(\mathcal{L}^{-A}), \qquad  e^{\mathcal{L}^ \alpha} \le  p_{\min} <  p_{\max} \le e^{\mathcal{L}^ \beta} < z.
\end{equation}
Each interval has essentially the same reciprocal sum. Its position among the other intervals matters when we select the active family.

Both the rough M\"obius estimate and the arithmetic decomposition require a bound for the contribution of large $z$-smooth divisors. At most $3x/a$ integers up to $3x$ have $a$ as a divisor, so we need the sum of $1/a$ beyond a threshold $T$.

\begin{lemma}[Smooth divisors]\label{vanu:smooth-tail}
We have that
\begin{align*}
  \sum_{P^ + (a) \le z} \frac1a = \prod_{p \le z} \Big(1 - \frac{1}{p} \Big)^{-1} \ll \log z.
\end{align*}
Moreover, uniformly for $T \ge 1$,
\begin{equation}  \label{vanu:smooth-tail-eq}
  \sum_{\substack{a > T \\P^ + (a) \le z}} \frac1a \ll (\log z) \exp \! \left( - \frac{\log T}{\log z} \right).
\end{equation}
\end{lemma}
\begin{proof}
The first bound follows from multiplicativity and Mertens theorem. For the tail, take $\eta = 1/ \log z$. When $a > T$ we have $1 \le (a/T)^ \eta$, so the tail \eqref{vanu:smooth-tail-eq} is at most $T^{-\eta} \prod_{p \le z}(1 - p^{-1 + \eta})^{-1}$. Then this product is $O(\log z)$ for sufficiently large $z$ by a Rankin calculation, see (2.9) in \cite[Proposition 2.2]{LWW}.
\end{proof}

\section{Prime-divisor patterns and exceptional integers}
\label{sec:dense-reconstruction}

The width $\gamma = \beta - \alpha$ will determine the exceptional density. Set
\begin{equation}  \label{gr:ambient-size}
  B = \left \lfloor \frac{\gamma \log \mathcal{L}}{\log2} \right \rfloor, \qquad \lambda_j = 2^{j - 1} \mathcal{L}^ \alpha, \qquad \sum_{j \le B} \lambda_j = (2^B - 1) \mathcal{L}^ \alpha < \mathcal{L}^ \beta.
\end{equation}
For sufficiently large $x$, one has $B \ge 1$ and $\tfrac12 \mathcal{L}^ \gamma < 2^B \le \mathcal{L}^ \gamma$. For $j \le B$, consider disjoint intervals of primes,
\begin{equation}  \label{gr:ambient-intervals}
  \mathcal{P}_j = \{p \in \mathbb P:e^{\lambda_j} < p \le e^{2 \lambda_j} \}.
\end{equation}

To detect when an interval contains a prime divisor of $n$, define
\begin{align*}
  \omega_j(n) = \sum_{p \in \mathcal{P}_j} \mathbf1_{p \mid n}, \qquad u_j(n) = \mathbf1_{\exists p \mid n : \, p \in \mathcal{P}_j} = \mathbf1_{\omega_j(n) \ge 1}.
\end{align*}

Fix $k \ge 1$. Call an integer $n$ \emph{typical} if $\sum_{j=1}^B u_j(n) \ge 3k$, and \emph{exceptional} otherwise; let $T,E$ be their indicators. Let $\mathcal{I}_{3k}$ consist of the $3k$-element subsets of $\{1,\ldots,B\}$. There are $|\mathcal{I}_{3k}|=\binom{B}{3k}\le B^{3k}$ patterns $\phi_J$, and we have an identity of disjoint indicators,
\begin{equation}  \label{gr:partition}
  T(n) = \sum_{J \in \mathcal{I}_{3k}}\phi_J(n).
\end{equation}
Here $\phi_J(n)=1$ means that $J$ gives the (unique) first $3k$ successful intervals of $n$.

For $J \in \mathcal{I}_{3k}$, enumerate its indices increasingly and put
\begin{align*}
  J = \{j_i:0 \le i < 3k\}, \qquad J^+ = \{j_{3i}:0 \le i < k\}, \qquad J^- = J \setminus J^+.
\end{align*}
Split $\phi_J = \phi_J^+\phi_J^-$ into active and passive factors,
\begin{align*}
  \phi_J^+(n) = \prod_{j \in J^+}u_j(n), \qquad \phi_J^-(n) = \prod_{j \in J^-}u_j(n)\prod_{\substack{1 \le j < \max J\\j \notin J}}(1-u_j(n)).
\end{align*}

For $j < j'$ in $J^+$, we have $j'-j \ge 3$, giving $\lambda_{j'}/\lambda_j \ge 8 > 19/3$.

For every prime $p$ in an active interval, the passive factor satisfies
\begin{equation}  \label{gr:spectator-invariance}
  |\phi_J^-(n)| \le 1, \qquad \phi_J^-(pn) = \phi_J^-(n).
\end{equation}

The exceptional-pair bound follows from sieving out two residue classes at each prime.

\begin{lemma}[Consecutive exceptions]\label{uniform:joint-residual}
For fixed $k,A > 0$ and $x \le y \le 2x$,
\begin{equation}  \label{uniform:joint-residual-bound}
  \sum_{x < n \le y}E(n)E(n + 1) \ll_k x \,2^{-2B}B^{6k - 2} + O_{A,k}(x \mathcal{L}^{-A}) \ll_{k,A}x \mathcal{L}^{-2 \gamma + o(1)} + O_{A,k}(x \mathcal{L}^{-A}).
\end{equation}
\end{lemma}
\begin{proof}
An exceptional integer hits at most $3k-1$ intervals. Set $s_0 = 6k-2$ and, for $S \subseteq \{1, \ldots,B \}$, define
\begin{align*}
  \mathcal{A}_S = \bigcup_{i \notin S} \mathcal{P}_i, \qquad \Psi_S(n) = \prod_{p \in \mathcal{A}_S} \mathbf1_{p \nmid n(n + 1)}.
\end{align*}
The union of the intervals hit by two exceptional neighbors has size at most $s_0$, so
\begin{equation}  \label{uniform:residual-majorization}
  E(n)E(n + 1) \le \sum_{\substack{S \subseteq \{1, \ldots,B \} \\|S| \le s_0}} \Psi_S(n).
\end{equation}
For squarefree $d\in\mathcal{A}_S$, the Chinese remainder theorem gives
\begin{equation}  \label{uniform:joint-divisor-count}
  \# \{x < n \le y:d \mid n(n + 1) \} = (y - x) \frac{2^{\omega(d)}}d + O(2^{\omega(d)}).
\end{equation}
The sieve with primes $\mathcal{A}_S$ has local density $g(p) = 2/p$. Then by Mertens' theorem,
\begin{align*}
  \prod_{\substack{v \le p < w \\p \in \mathcal{A}_S}}(1 - 2/p)^{-1} \ll \left(\frac{\log w}{\log v} \right)^2 \qquad(2 \le v \le w).
\end{align*}
Take $w = \exp(2 \mathcal{L}^ \beta)$ and $D = x^{1/3}$. All primes in $\mathcal{A}_S$ lie below $w$, and $\sigma = \log D/ \log w = \mathcal{L}^{1 - \beta}/6 \ge 19$ for large $x$. The fundamental lemma of sieve theory \cite[Theorem~6.9]{opera} gives
\begin{equation}  \label{uniform:joint-sieve-average}
  \sum_{x < n \le y} \Psi_S(n) = (y - x)V_S + O_A(x \mathcal{L}^{-A}), \qquad V_S = \prod_{p \in \mathcal{A}_S}(1 - 2/p),
\end{equation}
uniformly in $S$. Indeed, its error is
\begin{align*}
  \ll x e^{-\sigma} + \sum_{d \le D}2^{\omega(d)} \ll x e^{-\sigma} + D \log D \ll_A x \mathcal{L}^{-A}.
\end{align*}
Using \eqref{eq:interval-mass} and $\sum_{p > e^{\mathcal{L}^ \alpha}}p^{-2} \ll e^{-\mathcal{L}^ \alpha}$, we have
\begin{align*}
  V_S = 2^{-2(B - |S|)}(1 + o(1)) \ll_k2^{-2B} \qquad(|S| \le s_0).
\end{align*}
Sum \eqref{uniform:joint-sieve-average} over the $O_k(B^{s_0})$ choices in \eqref{uniform:residual-majorization} to complete the proof. Finally, note \eqref{gr:ambient-size} gives $2^{-2B}B^{s_0} = \mathcal{L}^{-2 \gamma + o(1)}$.
\end{proof}

The same argument, with the single forbidden residue $0$ at each prime and $|S| \le 3k-1$, gives
\begin{align*}
  \sum_{n \in \mathcal{J} \cap \mathbb Z}E(n) \ll_k x \,2^{-B}B^{3k - 1} + O_{A,k}(x \mathcal{L}^{-A}) \ll_k x \mathcal{L}^{-\gamma + o(1)} + O_{A,k}(x \mathcal{L}^{-A}).
\end{align*}
uniformly over intervals $\mathcal{J} \subseteq[x,3x]$.

\section{Arithmetic decomposition and Type~I estimates}
\label{sec:typeI}  \label{sec:arithmetic-reduction}

Here $n$ and $n+1$ play different roles. The smooth part of $n$ is captured by $\mathcal{U}$. Then a pattern $\phi$ captures the the first $3k$ successful intervals for $n + 1$. We decompose the factor at $n$ so that, in the Type~II sum, $\mathcal{U}$ is absorbed into the second coefficient. This leaves a (rough) M\"obius factor in the first coefficient for the major arc estimate. The Type~I sum reduces to M\"obius sums in arithmetic progressions.

Write every positive integer uniquely as
\begin{align*}
  n = s_n r_n, \qquad P^ + (s_n) \le z < P^-(r_n).
\end{align*}
Let $\mathcal{U}$ be any complex-valued function satisfying
\begin{equation}  \label{uniform:smooth-part-weight}
  | \mathcal{U}(n)| \le 1, \qquad \mathcal{U}(n) = \mathcal{U}(s_n).
\end{equation}
Truncate the smooth factor at $W$ by setting
\begin{equation}  \label{vanu:cutoffs}
  W = x^{1/100}, \qquad f_z(n) = \mu(n) \mathcal{U}(n) \mathbf1_{P^ + (n) \le z} \mathbf1_{n \le W}, \qquad F(n) = f_z(s_n)\mu(r_n).
\end{equation}

\begin{lemma}[Weighted approximation]
\label{vanu:approximation}  \label{uniform:first-weight}
For every fixed $A > 0$,
\begin{equation}  \label{vanu:approximation-eq}
  \sum_{n \le 3x}|F(n) - \mu(n) \mathcal{U}(n)| \ll_A x \mathcal{L}^{-A},
\end{equation}
uniformly in $\mathcal{U}$.
\end{lemma}
\begin{proof}
Since $(s_n,r_n) = 1$, we have $F(n) = \mathbf1_{s_n \le W}\mu(n)\mathcal{U}(n)$. Thus Lemma~\ref{vanu:smooth-tail} gives
\begin{align*}
  \sum_{n \le 3x}|F(n) - \mu(n)\mathcal{U}(n)| &\le \sum_{n \le 3x}\mathbf1_{s_n > W} \le 3x \sum_{\substack{u > W \\P^ + (u) \le z}} \frac1u \\
  &\ll x \mathcal{L}^{\zeta} \exp\!\left(-\frac1{100}\mathcal{L}^{1-\zeta}\right) \ll_A x\mathcal{L}^{-A}.
\end{align*}
\end{proof}

\begin{proposition}[Vaughan decomposition]
\label{prop:vau}  \label{vanu:decomposition}
Define
\begin{equation}  \label{vanu:coefficients}
  a(n) = \nu(n) \mathbf1_{n \le W^2}, \qquad b(n) = \nu(n) \mathbf1_{n > W^2}, \qquad A_q = (a*a*f_z)(q), \qquad \mathcal{B} = b*R_z*f_z.
\end{equation}
Then uniformly in $I \subseteq[x,3x]$,  $|G| \le 1$, for every fixed $A_0 > 0$,
\begin{equation}  \label{vanu:decomposition-eq}
  \sum_{n \in I \cap \mathbb Z} \mu(n) \mathcal{U}(n)G(n) = - S_{\rm I}(G;I) + S_{\rm II}(G;I) + O_{A_0}(x \mathcal{L}^{-A_0}),
\end{equation}
where
\begin{equation}  \label{vanu:type-I}
  S_{\rm I}(G;I) = \sum_{q \le W^5}A_q \sum_{qh \in I}R_z(h)G(qh)
\end{equation}
and
\begin{equation}  \label{vanu:type-II}
  S_{\rm II}(G;I) = \sum_{\substack{m,r > W^2\\mr \in I}} \mu(m)R_z(m)\mathcal{B}(r)G(mr).
\end{equation}
The coefficients satisfy
\begin{align*}
  |A_q| \le \tau(q), \qquad |\mathcal{B}(r)| \le \tau(r), \qquad A_q = 0 \ (q > x^{1/20}), \qquad \mathcal{B}(r) = 0 \ (r \le W^2).
\end{align*}
\end{proposition}
\begin{proof}
The support conditions give $(\nu*f_z)(n) = \mu(r_n)f_z(s_n) = F(n)$. The identity $\nu*R_z = \mathbf1_{\{1\}}$ follows at prime powers. Since $\nu = a + b$,
\begin{align*}
  \nu = (a + b)*(a + b)*R_z, \qquad a*a*R_z + a*b*R_z = a.
\end{align*}
Eliminating $a*b*R_z$ and convolving with $f_z$ gives
\begin{equation}  \label{vanu:exact-identity}
  F = 2a*f_z - A*R_z + b*\mathcal{B},
\end{equation}
where $A(q) = A_q$. The first term is supported on $n \le W^3 = x^{3/100} < x$. The other two give \eqref{vanu:type-I} and \eqref{vanu:type-II}, with the error bounded by Lemma~\ref{vanu:approximation}.

The smooth--rough decomposition yields
\begin{align}
  A_q & = \mu(s_q) \mathcal{U}(s_q) \mathbf1_{s_q \le W} \sum_{\substack{d_1d_2 = r_q \\d_1,d_2 \le W^2}} \mu(d_1) \mu(d_2), \label{vanu:coefficient-formulas} \\
  \mathcal{B}(t) & = \mu(s_t) \mathcal{U}(s_t) \mathbf1_{s_t \le W} \sum_{\substack{d \mid r_t \\d > W^2}} \mu(d). \notag
\end{align}
Plugging these formulae back in completes the proof.
\end{proof}

The decomposition gives $|\mathcal{B}(r)| \le \tau(r)$, whereas the bilinear estimate requires a pointwise bound by a power of $\mathcal{L}$. To discard the few coefficients above that threshold we use the divisor bounds
\begin{equation}  \label{vanu:divisor-moments}
  \sum_{n \le T} \tau_j(n) \ll_j T(\log(2T))^{j - 1}, \qquad \sum_{n \le T} \frac{\tau_j(n)}n \ll_j(\log(2T))^j, \qquad \tau_j(n)^2 \le \tau_{j^2}(n).
\end{equation}
The first and third bounds follow from \cite[Lemmas~1--2]{Mitkin}; the second follows by partial summation.

\begin{lemma}[Truncating the second coefficient]
\label{vanu:coefficient-truncation}
Given $C > 0$, then uniformly for $I \subseteq[x,3x]$ and $|G| \le 1$, replacing $\mathcal{B}(r)$ by $\mathcal{B}(r) \mathbf1_{|\mathcal{B}(r)| \le \mathcal{L}^C}$ changes $S_{\rm II}(G;I)$ by at most $O(x \mathcal{L}^{4 - C})$.
\end{lemma}
\begin{proof}
Since $|\mathcal{B}(r)| \le \tau(r)$, and $m$ takes at most $3x/r$ possible values, the change is bounded by
\begin{align*}
  3x \mathcal{L}^{-C} \sum_{r \le 3x} \frac{\tau(r)^2}{r} \le 3x \mathcal{L}^{-C} \sum_{r \le 3x} \frac{\tau_4(r)}r \ll x \mathcal{L}^{4 - C}.
\end{align*}
\end{proof}

For $Y \ge 2$, write
\begin{align*}
  E_\mu(Y;q) = \max_{\substack{a \;({\rm mod}\,q)\\t \le Y}} \bigg| \sum_{\substack{n \equiv a \;({\rm mod}\,q)\\n \le t }}\mu(n) \bigg|.
\end{align*}
\begin{lemma}[Weighted Bombieri--Vinogradov for $\mu$] \label{lem:mobius-BV}
For every $A > 0$, there exists $B(A) > 0$ such that
\begin{equation}  \label{rev:weighted-mobius-BV}
  \sum_{q \le Q} \tau_4(q)E_\mu(Y;q) \ll_A Y(\log Y)^{-A}.
\end{equation}
for all $Q \le Y^{1/2}(\log Y)^{-B(A)}$
\end{lemma}
\begin{proof}
This follows by \cite[Proposition 3]{MurVat18}, when combined with the divisor bound \ref{vanu:divisor-moments}.
\end{proof}

For disjoint subsets $\mathcal{U}, \mathcal{V} \subseteq \{1, \ldots,B \}$, set
\begin{align*}
  \phi_{\mathcal{U}, \mathcal{V}}(N) = \prod_{j \in \mathcal{U}}u_j(N) \prod_{j \in \mathcal{V}}(1 - u_j(N)), \qquad s = | \mathcal{U}| + | \mathcal{V}|.
\end{align*}
Fix an integer $R_0 \ge 1$ and indices $j_1 < \cdots < j_s$. Define
\begin{align*}
  Q_i(N) &= \begin{cases}
    \mathbf1_{\omega_{j_i}(N) \ge 1}, & j_i \in \mathcal{U}, \\[2pt]
    \mathbf1_{\omega_{j_i}(N) = 0}, & j_i \in \mathcal{V},
  \end{cases}\qquad T_i(N) = \begin{cases}
    \displaystyle \sum_{r = 1}^{R_0} (-1)^{r + 1}\binom{\omega_{j_i}(N)}{r}, & j_i \in \mathcal{U}, \\[2pt]
    \displaystyle \sum_{r = 0}^{R_0} (-1)^r\binom{\omega_{j_i}(N)}{r}, & j_i \in \mathcal{V}.
  \end{cases}
\end{align*}
Brun's sieve \cite[Chapter~6]{opera} gives $|Q_i(N) - T_i(N)| \le \binom{\omega_{j_i}(N)}{R_0 + 1}$. Since $\phi_{\mathcal{U}, \mathcal{V}}(N) = \prod_{i = 1}^s Q_i(N)$ and $0 \le Q_i(N) \le 1$, telescoping gives
\begin{align*}
  \Big| \phi_{\mathcal{U}, \mathcal{V}}(N) - \prod_{i = 1}^s T_i(N)\Big| &= \left|\sum_{i = 1}^s \bigl(Q_i(N) - T_i(N)\bigr)\prod_{l < i} T_l(N)\prod_{l > i} Q_l(N)\right|\\
  &\le \sum_{i = 1}^s |Q_i(N) - T_i(N)|\prod_{l < i} |T_l(N)| \le \sum_{i = 1}^s \binom{\omega_{j_i}(N)}{R_0 + 1}\prod_{l < i}\left(\sum_{r = 0}^{R_0}\binom{\omega_{j_l}(N)}{r}\right).
\end{align*}
Expand the main product and this nonnegative majorant into divisor sums with coefficients $c_d$ and $e_d$, respectively. Disjointness of the prime intervals gives $c_d \in \{0, \pm1 \}$. In the $i$th error term, the last nonempty interval is $j_i$ and contributes exactly $R_0 + 1$ prime factors; earlier intervals each contribute at most $R_0$. Thus the error terms have disjoint divisor supports, giving $e_d \in \{0,1 \}$ and
\begin{equation}  \label{gr:finite-expansion-eq}
  \left| \phi_{\mathcal{U}, \mathcal{V}}(N) - \sum_{d \mid N}c_d \right| \le \sum_{d \mid N}e_d.
\end{equation}
For $s = 0$, take $c_1 = 1$, all other $c_d = 0$, and $e_d = 0$. All supported divisors are odd and squarefree, with prime factors in the constrained intervals. By \eqref{gr:ambient-size},
\begin{equation}  \label{gr:divisor-support}
  d \le D_{R_0}: = \exp \bigl(2(R_0 + 1) \mathcal{L}^ \beta \bigr), \qquad \omega(d) \le sR_0 + 1.
\end{equation}
The primes in each interval have reciprocal sum at most one, so products of $r$ distinct primes from that interval have reciprocal sum at most $1/r!$. Thus the supported divisors in the main product and in the $i$th error term have reciprocal sums at most $e^s$ and $e^{i - 1}/(R_0 + 1)!$, respectively. Hence
\begin{equation}  \label{gr:divisor-masses}
  \sum_d \frac{|c_d|}{d} \le e^s, \qquad \sum_d \frac{e_d}{d} \le \frac{s e^s}{(R_0 + 1)!}.
\end{equation}
Since $s \le B < \log \mathcal{L}$, taking $R_0 = \lceil C \log \mathcal{L} \rceil$ for any fixed $C > 0$ gives $D_{R_0} = x^{o(1)}$, $e^s \le \mathcal{L}$, and, by Stirling's formula, $\sum_d e_d/d = O_A(\mathcal{L}^{-A})$ for every fixed $A$.

\begin{proposition}[Type~I cancellation]
\label{rev:mobius-typeI}  \label{uniform:typeI}
Let $G(n)=\mu(n+1)\phi(n+1)$ where $\phi = \phi_{\mathcal{U}, \mathcal{V}}$ is any prime-interval pattern defined above, and suppose that $a_q$ is supported on $q \le x^{1/20}$ with $|a_q| \le \tau(q)$. For every fixed $A > 0$, we have
\begin{equation}  \label{rev:mobius-typeI-general}
  S_{\rm I}(G;I)=\sum_{q \le x^{1/20}}a_q \sum_{qh \in I}R_z(h) \mu(qh + 1) \phi(qh + 1) \ll_A x \mathcal{L}^{-A}.
\end{equation}
The same estimate holds with $qh + 1$ replaced by $qh - 1$.
\end{proposition}
\begin{proof}
Choose $R = \lceil C_0 \log \mathcal{L} \rceil$, with $C_0 > e$ sufficiently large in terms of $A$. Finite inclusion--exclusion gives
\begin{equation}  \label{rev:sieve-truncation}
  R_z(h) = \sum_{\substack{\omega(r) \le R \\P^ + (r) \le z}} \mu(r) \mathbf{1}_{r \mid h} + O \left(\sum_{\substack{\omega(r) = R + 1 \\P^ + (r) \le z}}| \mu(r) | \mathbf{1}_{r \mid h} \right).
\end{equation}
Every divisor here is at most $z^{R + 1} = x^{o(1)}$. Apply \eqref{gr:finite-expansion-eq} with $R_0 = R$, using the coefficients $c_d,e_d$ above. Their support satisfies $d \le D_R = x^{o(1)}$.

If $r \mid h$ and $d \mid qh + 1$, then $(d,qr) = 1$. Hence writing $h = r \ell$ leaves one residue class for $\ell$ mod $d$, so
\begin{equation}  \label{rev:joint-divisibility-count}
  \# \{h:qh \in I, \ r \mid h, \ d \mid qh + 1 \} \le \frac{2x}{qrd} + 1.
\end{equation}
Insert \eqref{rev:sieve-truncation} first, bounding its error with $| \phi| \le 1$, and then expand $\phi$ in the remaining sum. Set $S_z = \sum_{p \le z}p^{-1}$. We have
\begin{align*}
  \sum_{q \le x^{1/20}} \frac{|a_q|}{q} \ll \mathcal{L}^2, \qquad \sum_{q \le x^{1/20}}|a_q| \ll x^{1/20} \mathcal{L}, \qquad \sum_r \frac1r \le \prod_{p \le z}(1 + p^{-1}) \ll \log z.
\end{align*}
Using \eqref{gr:divisor-masses} and \eqref{rev:joint-divisibility-count}, the total error is
\begin{equation}  \label{rev:pattern-truncation-error}
  \ll x \mathcal{L}^2 \left \{\frac{S_z^{R + 1}}{(R + 1)!} + (\log z) \frac{s e^s}{(R + 1)!} \right \} + x^{1/20 + o(1)}.
\end{equation}
The last term comes from the $ + 1$ in the progression count. For each $q$ there are only $x^{o(1)}$ choices of $(r,d)$. Since $\zeta < 1$, Mertens' estimate gives $S_z \le \log \mathcal{L}$ for sufficiently large $x$, and hence
\begin{align*}
  \frac{S_z^{R + 1}}{(R + 1)!} \le (e/C_0)^{R + 1}.
\end{align*}
The other factorial ratio has every fixed logarithmic saving by Stirling's formula and $s \le B < \log \mathcal{L}$. Taking $C_0$ sufficiently large makes \eqref{rev:pattern-truncation-error} $O_A(x \mathcal{L}^{-A})$.

\proofstep{Applying Bombieri--Vinogradov.}
Writing $h = r\ell$ and $N = qr\ell + 1$, the remaining main term is
\begin{align*}
  \sum_{\substack{q,r,d\\(d,qr)=1}}a_q\mu(r)c_d\sum_{\substack{N\in I+1\\N\equiv1\ ({\rm mod}\ qr)\\N\equiv0\ ({\rm mod}\ d)}}\mu(N),
\end{align*}
where $q,r,d$ range over their respective supports in the expansion. The Chinese remainder theorem makes the inner sum one residue class modulo $m=qrd$. For every fixed $B$ and sufficiently large $x$, these moduli satisfy
\begin{align*}
  m\le Q_*:=x^{1/20}z^R D_R=x^{1/20+o(1)}\le(4x)^{1/2}(\log(4x))^{-B}.
\end{align*}
Grouping triples by $qrd=m$ gives total absolute coefficient at most
\begin{align*}
  \sum_{qrd=m}|a_q\mu(r)c_d|\le\sum_{qrd=m}\tau(q)=\tau_4(m).
\end{align*}
Each interval sum is bounded by $2E_\mu(4x;m)$, so Lemma~\ref{lem:mobius-BV} bounds the main term by
\begin{align*}
  2\sum_{m \le Q_*}\tau_4(m)E_\mu(4x;m) \ll_A x\mathcal{L}^{-A}.
\end{align*}
Together with \eqref{rev:pattern-truncation-error}, this proves \eqref{rev:mobius-typeI-general}. For the mirrored sum with $qh-1$, take $N=qr\ell-1\in I-1$ and replace $N\equiv1\pmod{qr}$ by $N\equiv-1\pmod{qr}$; the same estimate applies.
\end{proof}

\section{The determinant bilinear form}  \label{rev:sec-typeII}  \label{sec:comparison}  \label{sec:determinant}

On the minor arcs, we do not exploit any M\"obius cancellation, and only use the variable size and support on rough integers. The estimate will be stated in Proposition~\ref{prop:min} and proved in Sections~\ref{sec:minor}--\ref{sec:minor-moments}.

Throughout, $\mathcal{L} = \log x$ and the prime intervals and roughness cutoff are those fixed in Section~\ref{sec:prelim}. We allow the dyadic lengths to satisfy
\begin{equation}  \label{eq:comparison-scales}
  M,N \ge x^ \delta, \qquad x/4 \le X: = MN \le 3x,
\end{equation}
where $\delta > 0$ is fixed. Thus $\log X \asymp \mathcal{L}$. Write $\mathcal{P}_i$ for the active intervals, $V_i = \sum_{p \in \mathcal{P}_i}1/p$, and $V = \prod_{i=1}^k V_i$. We shall group products of the chosen primes by their dyadic scale $Y$. The relevant scales satisfy
\begin{equation}  \label{eq:shift-scale-separation}
  V_i = \log2 + o(1), \qquad p_{\min} > e^{\mathcal{L}^ \alpha}, \qquad \mathcal{L}^ \alpha \ll \log Y \ll_k \mathcal{L}^ \beta, \qquad Y = o(z).
\end{equation}
In particular, a nonzero determinant of size $O(Y)$ cannot have a prime divisor above $z$.

\label{sec:cauchy}  \label{rev:indicator-typeii}
Fix an admissible family of $k$ active prime intervals. Set
\begin{align*}
  \mathcal{A} = \{p_1 \cdots p_k:p_i \in \mathcal{P}_i\}, \qquad \phi^+(v) = \prod_{i = 1}^k \mathbf1_{\omega_i(v) \ge 1}, \qquad \mathcal{P}_0 = \bigcup_{i=1}^k \mathcal{P}_i.
\end{align*}
Write $\phi = \phi^+\phi^-$, where for $v \ge 1$ and $p \in \mathcal{P}_0$, the passive factor $\phi^-$ satisfies
\begin{equation}  \label{rev:spectator-invariance}
  |\phi^-(v)| \le 1, \qquad \phi^-(pv) = \phi^-(v).
\end{equation}
\begin{lemma}[Ramar\'e identity]\label{lem:ramare-identity}
For all $v \ge 1$,
\begin{equation}  \label{rev:ramare-total-mass}
  \sum_{\substack{a \mid v \\a \in \mathcal{A}}} \frac1{\prod_i(\omega_i(v/a) + 1)} \le 1.
\end{equation}
If $p^2 \nmid v$ for all $p \in \mathcal{P}_0$, then
\begin{equation}  \label{rev:ramare-identity}
  \mu(v)\phi(v) = (-1)^k \sum_{\substack{a \mid v \\a \in \mathcal{A}}} \frac{\mu(v/a)\phi^-(v/a)}{\prod_i(\omega_i(v/a) + 1)}.
\end{equation}
\end{lemma}
\begin{proof}
If $\phi^+(v) = 0$, no $a \in \mathcal{A}$ divides $v$, so both claims hold. Otherwise, for $a = p_1 \cdots p_k \mid v$ with $p_i \in \mathcal{P}_i$, disjointness gives
\begin{align*}
  \omega_i(v/a) + 1 = \omega_i(v) + \mathbf1_{p_i^2 \mid v}.
\end{align*}
Choosing $a$ amounts to choosing one prime divisor from each interval. Hence
\begin{align*}
  \sum_{\substack{a \mid v \\a \in \mathcal{A}}} \frac1{\prod_i(\omega_i(v/a) + 1)} = \prod_{i=1}^k \sum_{\substack{p \mid v \\p \in \mathcal{P}_i}} \frac1{\omega_i(v) + \mathbf1_{p^2 \mid v}} \le \prod_{i=1}^k \frac{\omega_i(v)}{\omega_i(v)} = 1.
\end{align*}
If $p^2 \nmid v$ for all $p \in \mathcal{P}_0$, equality holds and $(a,v/a) = 1$. Thus multiplicativity and \eqref{rev:spectator-invariance} give
\begin{align*}
  \mu(v/a) = \frac{\mu(v)}{\mu(a)} = (-1)^k\mu(v), \qquad \phi^-(v/a) = \phi^-(v).
\end{align*}
Since $\phi^+(v) = 1$, we obtain
\begin{align*}
  \mu(v)\phi(v) = \mu(v)\phi^-(v) \sum_{\substack{a \mid v \\a \in \mathcal{A}}} \frac1{\prod_i(\omega_i(v/a) + 1)} = (-1)^k \sum_{\substack{a \mid v \\a \in \mathcal{A}}} \frac{\mu(v/a)\phi^-(v/a)}{\prod_i(\omega_i(v/a) + 1)}.
\end{align*}
\end{proof}

Consider the Type II sum
\begin{align}
  T = \sum_{m \sim M} \sum_{n \sim N} \alpha_m \beta_n \mu(mn + 1) \phi(mn + 1), \label{eq:Tdef}
\end{align}
The coefficients are supported on $m \sim M$, $n \sim N$, with $|\alpha_m|,|\beta_n| \le \mathcal{L}^C$ and $\alpha_m = 0$ unless $R_z(m) = 1$.

For $p \in \mathcal{P}_0$ and supported $m$, roughness gives $(m,p) = 1$, so $p^2 \mid mn + 1$ fixes $n \equiv -\overline m \pmod{p^2}$. Since $p^2 = x^{o(1)} \ll N$,
\begin{align}  \label{rev:indicator-square-error}
  \sum_{\substack{m \sim M,\ n \sim N\\\exists p \in \mathcal{P}_0:p^2 \mid mn + 1}} |\alpha_m\beta_n| &\ll X\mathcal{L}^{2C}\sum_{p \in \mathcal{P}_0}p^{-2} \ll X\mathcal{L}^{2C}e^{-\mathcal{L}^{\alpha}} \ll_A X\mathcal{L}^{-A}.
\end{align}
Set $c(q) = \mu(q)\phi^-(q)/\prod_i(\omega_i(q) + 1)$, and so $|c(q)| \le 1$. Since $\mu$ vanishes on squares, \eqref{rev:ramare-identity} gives the first equality; extending after extraction costs \eqref{rev:indicator-square-error} by \eqref{rev:ramare-total-mass}:
\begin{align}
  T &= (-1)^k \underset{p^2 \nmid mn + 1\ \forall p \in \mathcal{P}_0}{\sum_{m \sim M}\sum_{n \sim N}} \alpha_m\beta_n \sum_{\substack{a \mid mn + 1\\a \in \mathcal{A}}}c\!\left(\frac{mn + 1}{a}\right) \nonumber \\
  &= (-1)^k \sum_{\substack{a \in \mathcal{A},\ m,n\\a \mid mn + 1}} \alpha_m\beta_n c\!\left(\frac{mn + 1}{a}\right) + O_A(X\mathcal{L}^{-A}). \label{eq:T2}
\end{align}
Choose a fixed nonnegative smooth function $\eta$, supported on $[1,4]$, whose rescalings $\eta_Y(y) = \eta(y/Y)$ satisfy
\begin{align*}
  \sum_{Y \in 2^{\mathbb Z}}\eta_Y(a) = 1 \qquad(a > 0).
\end{align*}
Only $O_k(\mathcal{L})$ scales $Y$ occur, with $\mathcal{L}^{\alpha} \ll \log Y \ll_k \mathcal{L}^{\beta}$. Define
\begin{align*}
  S_Y(q) = \sum_{\substack{a \in \mathcal{A},\ m,n\\mn + 1 = aq}} \eta_Y(a)\alpha_m\beta_n, \qquad T_Y = (-1)^k \sum_{q \ge 1}c(q)S_Y(q).
\end{align*}
Inserting the partition into \eqref{eq:T2} and grouping by $q$ gives
\begin{align}  \label{eq:T3}
  T = \sum_Y T_Y + O_A(X\mathcal{L}^{-A}).
\end{align}
Since $a \ge Y$ and $mn \le 4X$ on the support, $S_Y(q) = 0$ for $q > (4X + 1)/Y$. Thus Cauchy--Schwarz gives
\begin{equation}  \label{eq:CS}
  |T_Y|^2 \le \sum_{q \le (4X + 1)/Y}|c(q)|^2 \cdot \sum_{q \ge 1}|S_Y(q)|^2 \ll \frac XY \sum_{q \ge 1}|S_Y(q)|^2.
\end{equation}
Expanding the square gives
\begin{align}  \label{eq:013}
  S_Y := \sum_{q \ge 1}|S_Y(q)|^2 = \sum_{\substack{a,b \in \mathcal{A},\ m,n,r,s\\(mn + 1)/a = (rs + 1)/b \in \mathbb N}} \eta_Y(a)\eta_Y(b)\alpha_m\overline{\alpha_r}\beta_n\overline{\beta_s}.
\end{align}

We first remove pairs with $(a,b) > 1$ to replace the cofactor equality by a determinant equation.

\begin{lemma}[Common chosen primes]\label{lem:shared}
Under \eqref{eq:comparison-scales}, the total absolute contribution to $S_Y$ in \eqref{eq:013} of pairs $a,b$ sharing an prime $p\in\mathcal{P}_0$ is $O_{A,k}(XY \mathcal{L}^{-A})$ for every fixed $A > 0$.
\end{lemma}
\begin{proof}
Denote this contribution by $E_Y$. We modify the divisor estimate from \cite[Lemma~4.1]{LWW}.

Write $\tau_j$ for the $j$-fold divisor function. Choosing a largest factor in each ordered four-factor decomposition gives
\begin{align*}
  \tau(u)^2 \le \tau_4(u) \le 4 \sum_{\substack{d \mid u \\d \le u^{3/4}}} \tau_3(d).
\end{align*}
Also $d \mid \ell q - 1$ implies $(d, \ell) = 1$ and $q \equiv \overline \ell \pmod d$. Thus we have, uniformly for $\ell \asymp Y = x^{o(1)}$,
\begin{equation}  \label{eq:014}
  \sum_{\substack{q \ge 1 \\X < \ell q - 1 \le 4X}} \tau(\ell q - 1)^2 \ll \sum_{d \le (4X)^{3/4}} \tau_3(d) \left(\frac{X}{Yd} + 1 \right) \ll \frac XY \mathcal{L}^3 + X^{3/4} \mathcal{L}^2 \ll \frac XY \mathcal{L}^3.
\end{equation}

For fixed $a,b \in \mathcal{A} \cap[Y,4Y]$, we count the factorizations of $aq - 1$ and $bq - 1$, apply Cauchy--Schwarz and \eqref{eq:014} to give
\begin{align*}
  \sum_{\substack{q \ge 1, \ m,n,r,s \\mn + 1 = aq, \ rs + 1 = bq}} &| \alpha_m \alpha_r \beta_n \beta_s| \le \mathcal{L}^{4C} \sum_{\substack{ X < aq - 1 \le 4X \\X < bq - 1 \le 4X}} \tau(aq - 1) \tau(bq - 1) \\
  & \le \mathcal{L}^{4C} \left(\sum_{X < aq - 1 \le 4X} \tau(aq - 1)^2 \right)^{1/2} \left(\sum_{X < bq - 1 \le 4X} \tau(bq - 1)^2 \right)^{1/2} \ll \frac XY \mathcal{L}^{4C + 3}.
\end{align*}

The number of pairs sharing a prime satisfies
\begin{equation}  \label{eq:abcoll}
  \begin{aligned}
  \# \{a,b \in \mathcal{A} \cap[Y,4Y]:(a,b) > 1 \} &\le \sum_{p \in \mathcal{P}_0} \left(\sum_{p \mid a \in \mathcal{A} \cap[Y,4Y]}1 \right)^2 \le \sum_{p \in \mathcal{P}_0} \left \lfloor \frac{4Y}{p} \right \rfloor^2\\
  &\le 16Y^2 \sum_{p \in \mathcal{P}_0} \frac1{p^2} \ll Y^2e^{-\mathcal{L}^ \alpha}.
  \end{aligned}
\end{equation}
Since $| \eta_Y| \ll 1$ and $V^{-2} \ll_k1$, multiplying these estimates gives
\begin{align*}
  E_Y \ll_k \frac XY \mathcal{L}^{4C + 3} \cdot Y^2e^{-\mathcal{L}^ \alpha} \ll_{A,k}XY \mathcal{L}^{-A}.
\end{align*}
The diagonal $a = b$ is included.
\end{proof}

Note the equation $b(mn + 1) = a(rs + 1)$ implies $a \mid b(mn + 1)$, and so $a \mid mn + 1$ for $(a,b) = 1$. Similarly, $b \mid rs + 1$ for $(a,b) = 1$. Hence we have $bmn - ars = a - b$ if and only if
\begin{equation}  \label{eq:det2}
  \frac{mn + 1}{a} = \frac{rs + 1}{b} \in \mathbb N.
\end{equation}
Thus Lemma~\ref{lem:shared} yields
\begin{align*}
  S_Y = \sum_{\substack{a,b \in \mathcal{A}, \ (a,b) = 1 \\m,n,r,s \\bmn - ars = a - b}} \eta_Y(a) \eta_Y(b) \alpha_m \overline{\alpha_r} \beta_n \overline{\beta_s} \quad + \ O_{A,k}(XY \mathcal{L}^{-A}).
\end{align*}
Define the quantity $t=bmn-ars$. To detect $t = a - b$, we apply the Fourier identity,
\begin{align*}
  \mathbf1_{t = a - b} = \int_{\mathbb T} e \bigl(\theta(t + b - a) \bigr) \,d \theta.
\end{align*}

Since $a,b \in[Y,4Y]$, a solution has $|t| \le 3Y$. Thus we may insert an even real cutoff $\psi_1 \in C_c^ \infty(\mathbb R)$ equal to one on $[ - 4,4]$. Splitting into major and minor arcs,
\begin{equation}  \label{eq:arcs-kernels}
  \begin{aligned}
  \mathfrak M& = \bigcup_{\substack{1 \le q \le \mathcal{L}^{A_0} \\(j,q) = 1}} \{\theta \in \mathbb T: \| \theta - j/q \| \le 2 \mathcal{L}^{A_0}/Y \}, \qquad \text{and} \qquad \mathfrak m& = \mathbb T \setminus \mathfrak M, \\
  H^{\rm maj}(t;a,b)& = \psi_1(t/Y) \int_{\mathfrak M}e \bigl(\theta(t + b - a) \bigr) \,d \theta, \\
  H^{\rm min}(t;a,b)& = \psi_1(t/Y) \int_{\mathfrak m}e \bigl(\theta(t + b - a) \bigr) \,d \theta.
  \end{aligned}
\end{equation}
Thus
\begin{equation}  \label{eq:uniform-kernel-size}
  | \mathfrak M| \ll \mathcal{L}^{3A_0}/Y, \qquad | \mathfrak m| \le 1, \qquad |H^{\rm maj}(t;a,b)| \ll \mathcal{L}^{3A_0}/Y, \qquad |H^{\rm min}(t;a,b)| \ll 1.
\end{equation}
Both kernels vanish if $|t| > C_\psi Y$, and the two arc sets and the cutoff are invariant under reflection. For integer $t$ on the support of $\eta_Y(a) \eta_Y(b)$, we have
\begin{align*}
  \mathbf1_{t = a - b} = H^{\rm maj}(t;a,b) + H^{\rm min}(t;a,b).
\end{align*}
Indeed, when $t = a - b$ we have $|t| \le 3Y$, so $\psi_1(t/Y) = 1$.

Define the normalized bilinear form $Q_Y = V^{-2}S_Y$, and
\begin{align}
  Q_Y^{\rm maj} &= V^{-2} \sum_{\substack{a,b \in \mathcal{A} \\m,n,r,s}} \eta_Y(a) \eta_Y(b) \alpha_m \overline{\alpha_r} \beta_n \overline{\beta_s}H^{\rm maj}(bmn - ars;a,b), \label{eq:Qmaj}\\
  Q_Y^{\rm min} &= Q_Y - Q_Y^{\rm maj}.  \label{eq:Qmin}
\end{align}
For $(a,b) = 1$, Fourier orthogonality identifies the contribution to $Q_Y^{\rm min}$ with the sum using $H^{\rm min}$. For $(a,b) > 1$ the definition keeps the difference between the original cofactor sum and the major arc sum. Lemma~\ref{lem:shared} bounds the former contribution; \eqref{eq:abcoll} and the fixed-pair bound used in \eqref{eq:032} bound the latter.

\section{Completion of the proof of Theorem \ref{thm:chowla}}  \label{sec:final}

We complete the proof assuming the major and minor arc estimates 
\begin{align*}
  |Q_Y^{\rm maj}|, \ |Q_Y^{\rm min}| \ll_A XY \mathcal{L}^{-A}
\end{align*}
in Propositions~\ref{prop:maj} and~\ref{prop:min}. 

First we apply Perron's formula to separate the variables from the interval condition $mn \in I$

\begin{proposition}[Fixed-prime-interval Type II estimate]
\label{rev:uniform-typeii}  \label{uniform:unit-fixed-prime-interval}
Fix the parameters in \eqref{par:all} and $\delta,C,C_0,C_w,A > 0$. For sufficiently large fixed $k$, let the active prime intervals be admissible and suppose
\begin{align*}
  x/4 \le X = MN \le 3x, \qquad M,N \ge x^ \delta.
\end{align*}
Take coefficients on arbitrary subintervals of $(M,2M]$ and $(N,2N]$ with
\begin{align*}
  \alpha_m = \mu(m)R_z(m)m^{iv}w(m), \quad |v| \le x^{C_0}, \qquad | \beta_n| \le \mathcal{L}^C, \qquad \|w \|_\infty + \operatorname{Var}(w) \le C_w.
\end{align*}
For every passive factor $\phi^-$ satisfying \eqref{rev:spectator-invariance} and every summation interval $I$,
\begin{equation}  \label{rev:indicator-bilinear}
  \sum_{\substack{m,n \\mn \in I}} \alpha_m \beta_n \mu(mn + 1) \phi(mn + 1) \ll_A X \mathcal{L}^{-A}.
\end{equation}
\end{proposition}
\begin{proof}
Choose $A_0$, then $k$, in Propositions~\ref{prop:maj} and~\ref{prop:min}, replacing $(A,C_0,C_w)$ by $(2A + 20,\max(C_0,3) + 1,2C_w)$. Then $Q_Y \le |Q_Y^{\rm maj}| + |Q_Y^{\rm min}| \ll XY \mathcal{L}^{-2A - 20}$.

Write $T(\alpha,\beta)$ for the sum in \eqref{eq:Tdef}. Then by \eqref{eq:T3} and \eqref{eq:CS}, we have
\begin{align*}
  T(\alpha,\beta) &= (-1)^k \sum_Y \sum_q \frac{\mu(q)\phi^-(q)}{\prod_i(\omega_i(q) + 1)}S_Y(q) + O(X \mathcal{L}^{-A - 10}),\\
  |T(\alpha,\beta)| &\ll_k \sum_Y \left(\frac{X}{Y}Q_Y\right)^{1/2} + X \mathcal{L}^{-A - 10} \ll X \mathcal{L}^{-A - 9}.
\end{align*}
For $\sigma = 1/\mathcal{L}$, $s = \sigma + it$, and $|t| \le x^3$, set
\begin{align*}
  \alpha_m^{(s)} &:= \alpha_m m^{-s} = \mu(m)R_z(m)m^{i(v - t)}w_\sigma(m), & w_\sigma(u) &:= u^{-\sigma}w(u),\\
  \beta_n^{(s)} &:= \beta_n n^{-s}, & |\beta_n^{(s)}| &\le \mathcal{L}^C,\\
  |v - t| &\le x^{\max(C_0,3) + 1}, & \|w_\sigma\|_\infty + \operatorname{Var}(w_\sigma) &\le 2C_w.
\end{align*}
Choose half-integers $u_1,u_2 \asymp x$ such that $mn \in I$ iff $u_1 < mn < u_2$. Hence truncated Perron gives
\begin{align*}
  \sum_{\substack{m,n\\mn \in I}} \alpha_m \beta_n \mu(mn + 1) & \phi(mn + 1) = \frac{1}{2\pi i}\int_{\sigma - ix^3}^{\sigma + ix^3} \frac{u_2^s - u_1^s}{s}\,T(\alpha^{(s)},\beta^{(s)})\,ds + O(\mathcal{L}^C/x)\\
  &\ll X \mathcal{L}^{-A - 9}\int_{-x^3}^{x^3}\frac{dt}{|\sigma + it|} + \mathcal{L}^C/x \ll X \mathcal{L}^{-A}.
\end{align*}
\end{proof}

In the final sum the typical integer can be either member of a consecutive pair. For the preceding integer, replace $mn + 1$ by $mn - 1$ in $\mu$ and $\phi$.

\begin{corollary}[Vaughan Type II sum]\label{cor:vaughan-typeii}
Let $W$ be as in \eqref{vanu:cutoffs}. For fixed $C,A > 0$, every interval $I \subset[x,3x]$, and every coefficient $b(n)$ supported on $n > W^2$ with $|b(n)| \le \mathcal{L}^C$, we have uniformly
\begin{equation}  \label{rev:indicator-vaughan-II}
  \sum_{\substack{m,n > W^2 \\mn \in I}} \mu(m)R_z(m)b(n) \mu(mn + 1) \phi(mn + 1) \ll_A x \mathcal{L}^{-A}.
\end{equation}
The same bound also holds for the reflected form with $mn - 1$.
\end{corollary}
\begin{proof}
Divide $m,n > W^2$ into $O(\mathcal{L}^2)$ dyadic intervals. Each contributing box has $x/4 \le MN \le 3x$ and $M,N \ge W^2/2 \ge W = x^{1/100}$ for large $x$. Apply \eqref{rev:indicator-bilinear} with $A + 3$ and sum over intervals.
\end{proof}

We now combine the Type~I and Type~II estimates for a fixed pattern.

\begin{proposition}[Patterns]
\label{prop:completed-pattern}
Fix $A > 0$ and the exponents in \eqref{par:all}. For a sufficiently large fixed $k$, every pattern $\phi_J$ with $J \in \mathcal{I}_{3k}$ satisfies
\begin{equation}  \label{eq:completed-pattern}
  \sum_{n \in I \cap \mathbb Z} \mu(n) \mathcal{U}(s_n) \mu(n + 1) \phi_J(n + 1) \ll x \mathcal{L}^{-A}
\end{equation}
for every interval $I \subseteq[x,3x]$ and every weight $| \mathcal{U}| \le 1$. The bound is uniform in the interval, weight, active and passive intervals, with all prime intervals and cutoffs defined from the original base $x$. The same bound holds with $n + 1$ replaced by $n - 1$.
\end{proposition}
\begin{proof}
Set $G(n) = \mu(n + 1)\phi_J(n + 1)$ and apply Proposition~\ref{prop:vau} to the weight $n \mapsto \mathcal{U}(s_n)$ with precision $A + 4$:
\begin{align*}
  S := \sum_{n \in I \cap \mathbb Z}\mu(n)\mathcal{U}(s_n)\mu(n + 1)\phi_J(n + 1) = -S_{\rm I}(G;I) + S_{\rm II}(G;I) + O(x\mathcal{L}^{-A-4}).
\end{align*}
Since $|A_q| \le \tau(q)$ and $A_q = 0$ for $q > x^{1/20}$, Proposition~\ref{rev:mobius-typeI} gives $S_{\rm I}(G;I) \ll x\mathcal{L}^{-A-4}$.

Take $C > A + 12$ and set $\mathcal{B}_C(r) = \mathcal{B}(r)\mathbf1_{|\mathcal{B}(r)| \le \mathcal{L}^C}$. Then $|\mathcal{B}_C(r)| \le \mathcal{L}^C$, $\mathcal{B}_C(r) = 0$ for $r \le W^2$, and Lemma~\ref{vanu:coefficient-truncation} bounds the truncation error by $O(x\mathcal{L}^{4-C})$.

Write $\phi_J = \phi^+\phi^-$. The active scales satisfy
\begin{align*}
  \frac{\lambda_{j_{3i+3}}}{\lambda_{j_{3i}}} = 2^{j_{3i+3}-j_{3i}} \ge 8 > 19/3 \qquad(0 \le i < k-1),
\end{align*}
and \eqref{gr:spectator-invariance} gives $|\phi^-| \le 1$ and $\phi^-(pv) = \phi^-(v)$ for every active prime $p$. With $C$ fixed, choose $k$ sufficiently large and apply \eqref{rev:indicator-vaughan-II} at precision $A + 4$:
\begin{align*}
  S_{\rm II}(G;I) &= \sum_{\substack{m,r > W^2\\mr \in I}}\mu(m)R_z(m)\mathcal{B}_C(r)\mu(mr + 1)\phi^+(mr + 1)\phi^-(mr + 1) + O(x\mathcal{L}^{4-C}) \\
  &\ll x\mathcal{L}^{-A-4} + x\mathcal{L}^{4-C}.
\end{align*}
Hence $|S| \ll x\mathcal{L}^{-A-4} + x\mathcal{L}^{4-C} \ll x\mathcal{L}^{-A}$. For the preceding integer, take $G(n) = \mu(n - 1)\phi_J(n - 1)$ and use the reflected Type~I and Type~II estimates.
\end{proof}

\begin{corollary}[Typical-integer estimate]
\label{uniform:marked-pair}
For every fixed $A > 0$, a sufficiently large fixed $k$ gives
\begin{equation}  \label{uniform:marked-pair-eq}
  \sum_{n \in I \cap \mathbb Z} \mu(n) \mathcal{U}(s_n) \mu(n + 1)T(n + 1) \ll x \mathcal{L}^{-A},
\end{equation}
uniformly for $I \subseteq[x,3x]$ and $| \mathcal{U}| \le 1$. The same bound holds with $n + 1$ replaced by $n - 1$.
\end{corollary}
\begin{proof}
Apply Proposition~\ref{prop:completed-pattern} with $A + 2$, and sum over patterns,
\begin{equation}  \label{uniform:T-def}
  T(N) = \sum_{J \in \mathcal{I}_{3k}} \phi_J(N).
\end{equation}
For this fixed $k$, one has $B^{3k} \le \mathcal{L}$ for sufficiently large $x$.
\end{proof}

Now fix $0 < \eps < 1$ and take
\begin{equation}  \label{eq:final-parameters}
  \alpha = \frac{\eps}{8}, \qquad \beta = \frac12 - \frac{\eps}{8}, \qquad \jmath = \frac{\eps}{32}, \qquad \varrho = \frac12, \qquad \zeta = \frac{1 + \beta}{2}.
\end{equation}
These satisfy \eqref{par:all}; in particular,
\begin{align*}
  \beta + \jmath = \frac12 - \frac{3 \eps}{32} < \varrho < \frac12 + \frac{3 \eps}{32} = 1 - \beta - \jmath.
\end{align*}
Writing $\gamma = \beta - \alpha$, we have
\begin{equation}  \label{eq:final-width}
  2 \gamma = 1 - \frac{\eps}{2}.
\end{equation}
Choose $k$ from Corollary~\ref{uniform:marked-pair} with $A = 10$. All these choices remain fixed as $x$ tends to infinity.

Every available prime is below $z$, so
\begin{equation}  \label{uniform:smooth-residual-weight}
  E(n) = E(s_n), \qquad T(n) = T(s_n).
\end{equation}
Set $\mu_T = \mu T$ and $\mu_E = \mu E$. Once $n$ is typical the reflected estimate applies whatever happens at $n + 1$. Accordingly, we split according to $n + 1$ only when $n$ is exceptional. This gives the three terms
\begin{equation}  \label{uniform:two-residual-identity}
  \mu(n) \mu(n + 1) = \mu_T(n) \mu(n + 1) + \mu_E(n) \mu_T(n + 1) + \mu_E(n) \mu_E(n + 1).
\end{equation}
We sum this identity over $x < n \le y$, where $x \le y \le 2x$.

If $n$ is exceptional and $n + 1$ is typical, then apply Corollary~\ref{uniform:marked-pair} with the weight $\mathcal{U} = E$, permitted by \eqref{uniform:smooth-residual-weight}. It gives
\begin{align*}
  \sum_{x < n \le y} \mu_E(n) \mu_T(n + 1) \ll x \mathcal{L}^{-10}.
\end{align*}
If $n$ is typical, then writing $v = n + 1$ gives the reflected sum
\begin{align*}
  \sum_{x < n \le y} \mu_T(n) \mu(n + 1) = \sum_{x + 1 < v \le y + 1} \mu(v) \mu(v - 1)T(v - 1).
\end{align*}
Apply the reflected form of the same corollary with weight $1$, using the prime intervals defined from the original base $x$. This contribution is also $O(x \mathcal{L}^{-10})$.

If both $n,n+1$ are exceptional,
Lemma~\ref{uniform:joint-residual} gives
\begin{align*}
  \left| \sum_{x < n \le y} \mu_E(n) \mu_E(n + 1) \right| \le \sum_{x < n \le y}E(n)E(n + 1) & \ll_k x \,2^{-2B}B^{6k - 2} + O(x \mathcal{L}^{-10}) \\
  & \ll_k x \mathcal{L}^{-2 \gamma}(\log \mathcal{L})^{6k - 2} + O(x \mathcal{L}^{-10}).
\end{align*}
The choice of the remaining intervals costs only $(\log \mathcal{L})^{6k - 2}$, which is $O_{k, \eps}(\mathcal{L}^{\eps/4})$ since $k$ is fixed. By \eqref{eq:final-width} the exceptional pairs thus contribute $O_\eps(x \mathcal{L}^{-1 + 3 \eps/4})$. The two sums with a typical integer are $O(x\mathcal{L}^{-10})$, giving \eqref{eq:dyadic-goal}.

\begin{proof}[Proof of Theorem~\ref{thm:chowla}]
Partition $(\sqrt{X},X]$ into dyadic intervals $(x,2x]$ with $x \ge \sqrt X/2$. On these intervals $\log x \asymp \log X$, and their left endpoints sum to $O(X)$. Summing \eqref{eq:dyadic-goal} gives
\begin{align*}
  \left| \sum_{n \le X} \mu(n) \mu(n + 1) \right| \ll_\eps X(\log X)^{-1 + \eps} + \sqrt X \ll_\eps X(\log X)^{-1 + \eps}.
\end{align*}
\end{proof}

\section{The major arcs}  \label{sec:major}

On the major arcs we exploit cancellation in the rough M\"obius factor from the first coefficient. We adapt the proof of \cite[Proposition 5.3]{LWW} to allow unrestricted support for $\beta_n$ and the full range of prime intervals in \eqref{rev:admissible-intervals}. The desired cancellation ultimately derives from the Vinogradov--Korobov zero-free region for Dirichlet $L$-functions.

We shall use a M\"obius cancellation estimate of L\"u--Wang--Wu and OpenAI, itself a $z$-rough adaptation of Matom\"aki--Shao--Tao--Ter\"av\"ainen \cite[Lemma~3.9(ii)]{MSTT}.

\begin{proposition}[M\"obius--character cancellation]\label{prop:mel}
Let $\delta,B,C > 0$ be fixed. Suppose that $x^ \delta \le D \le 3x$ and $I \subseteq[D,2D]$ is an interval. Let $\chi$ be a Dirichlet character mod $q$ with $1 \le q \le \mathcal{L}^B$. Then, for any weight $w$ of bounded total variation, we have for every $A > 0$,
\begin{equation}  \label{eq:mellin-var}
  \sup_{|t| \le x^C} \left| \sum_{d \in I} \mu(d)R_z(d) \chi(d)d^{it}w(d) \right| \ll D \mathcal{L}^{-A} \bigl(\|w \|_\infty + \operatorname{Var}_I(w) \bigr).
\end{equation}
\end{proposition}
\begin{proof}
This is \cite[Proposition 2.2]{LWW} and \cite[(2.11)]{LWW}, with $\zeta$ instead of $0.24$ and $D\le 3x$ instead of $D\le 2x$. Both changes follow from their proof.
\end{proof}

Keep $x/4 \le X = MN \le 3x$, $M,N \ge x^ \delta$, and the notation of Section~\ref{sec:determinant}. Since each product $a \in \mathcal{A}$ has exactly one prime from each interval, the reciprocal sum over $\mathcal{A}$ factors by interval, giving
\begin{equation}  \label{eq:aharm}
  V^{-1} \sum_{a \in \mathcal{A}} \frac1a = 1.
\end{equation}

\begin{proposition}[Major arc estimate]\label{prop:maj}
Let $\alpha_m = \mu(m)R_z(m)m^{iv}w(m)$, $|v| \le x^{C_0}$, and $| \beta_n| \le \mathcal{L}^C$, with $\|w \|_\infty + \operatorname{Var}(w) \le C_w$. Then for $Q_Y^{\rm maj}$ in \eqref{eq:Qmaj}, we have
\begin{equation}  \label{eq:majbd}
  |Q_Y^{\rm maj}| \ll_A XY \mathcal{L}^{-A}.
\end{equation}
\end{proposition}
\begin{proof}
Write $\theta=j/q+\omega/Y$ on a major arc, with $q\le\mathcal{L}^{A_0}$ and $|\omega|\le2\mathcal{L}^{A_0}$. Since $b,m$ are units modulo $q$, character orthogonality gives
\begin{equation}  \label{maj:rational-expansion}
  \begin{aligned}
  e\!\left(\frac{jb(mn+1)}q\right)
    &=\sum_{\chi,\psi\bmod q}\Omega_\psi(1)(\chi\psi)(b)\chi(m)\Omega_\chi(n),\qquad
  \Omega_\chi(n)
    =\frac1{\varphi(q)}\sideset{}{^*}\sum_{u\bmod q}\overline{\chi(u)}e(jun/q).
  \end{aligned}
\end{equation}
Here $|\Omega_\chi(n)|\le1$ for every $n$, so $\beta_n\Omega_\chi(n)$ still has size at most $\mathcal{L}^C$. Applying \eqref{maj:rational-expansion} also at the conjugate endpoint gives at most $\varphi(q)^4$ terms.

By the Mellin reduction in \cite[proof of Proposition~5.3, (5.13)--(5.16)]{LWW}, it suffices at each endpoint to show that
\begin{align}\label{eq:major-arc-energy}
  \int_{-T}^T |P_\xi(t)M_\chi(t)N_\Omega(t)|^2\,dt
  \ll_H \mathcal{L}^{-H}
\end{align}
for every fixed $H,B>0$, where $T=2X\mathcal{L}^B$. Here $P_\xi$ comes from one prime interval, and
\begin{align*}
  P_\xi(t) &= \frac1P\sum_{P\le p<2P} v_p p^{i(t+\xi)},
  \qquad |v_p|\le1,\quad \xi\in\mathbb R,\quad
  \mathcal{L}^\alpha\ll\log P\ll\mathcal{L}^\beta,\\
  M_\chi(t) &= \frac1M\sum_m \alpha_m\chi(m)m^{it}, \qquad
  N_\Omega(t) = \frac1N\sum_n \beta_n\Omega_\chi(n)n^{it}.
\end{align*}
The remaining prime factors are bounded by their normalized reciprocal sums as in \eqref{eq:aharm}, and our smooth cutoff $\psi_1$ gives the same Fourier decay.

Let $\mathcal I$ be the unit intervals with integer endpoints meeting $\{|t|\le T:|P_\xi(t)|>\mathcal{L}^{-A_s}\}$. Take $\ell=\lfloor\log(4T)/\log(2P)\rfloor$. The squared coefficient norm of $P_\xi^\ell$ is at most $\ell!(2/P)^\ell$ by unique factorization. Choosing one point above the threshold per interval and splitting the intervals into three classes modulo $3$, the separated-point mean square \cite[Lemma~2.2, (2.3)]{PD} gives $|\mathcal I|\ll T\mathcal{L}^{O(1)}\ell!(2\mathcal{L}^{2A_s}/P)^\ell$. Hence
\begin{equation}  \label{maj:sparse-bounds}
  \log(2+|\mathcal I|)\ll\frac{\mathcal{L}\log\mathcal{L}}{\log P}+\log P
  \ll\mathcal{L}^{1-\alpha}\log\mathcal{L}+\mathcal{L}^\beta
  =o\!\left(\frac{\mathcal{L}}{(\log\mathcal{L})^2}\right).
\end{equation}
Thus the estimate in \cite[proof of Lemma~5.2, (5.9)--(5.11)]{LWW} gives $\sum_{I\in\mathcal I}\sup_{t\in I}|N_\Omega(t)|^2\ll\mathcal{L}^{2C+1}$.

Where $|P_\xi|\le\mathcal{L}^{-A_s}$, the product mean square \cite[(5.2)]{LWW} bounds the integral there by $O(\mathcal{L}^{-2A_s+2C+3+B})$. On the complementary set, $|P_\xi|\ll1$ and Proposition~\ref{prop:mel} gives $|M_\chi(t)|\ll_{A'}\mathcal{L}^{-A'}$, since $|t+v|$ is bounded by a fixed power of $x$. The last sum over $\mathcal I$ bounds this contribution by $O(\mathcal{L}^{-2A'+2C+1})$. Choosing $A_s,A'$ sufficiently large gives \eqref{eq:major-arc-energy} uniformly in $\xi$, and hence \eqref{eq:majbd}.
\end{proof}

\section{Prime lists and operators}  \label{sec:minor}

\begin{proposition}[Minor arc estimate]\label{prop:min}
Denote $Q_Y^{\rm min}$ in \eqref{eq:Qmin}. Then we have
\begin{equation}  \label{eq:minbd}
  |Q_Y^{\rm min}| \ll_A XY \mathcal{L}^{-A}.
\end{equation}
\end{proposition}

We adapt the determinant construction of \cite[Sections~4.2--4.7]{LWW} and \cite[Sections~3--4]{PD}, and the treatment of prime lifetimes in \cite[Sections~1.5 and~5]{HR2021} and \cite[Proposition~5.3]{Pilatte}.

\begin{lemma}[Exceptional endpoints]\label{lem:exceptional-endpoints}
(i) If $(a,b)=1$, then $bmn-ars=a-b$ implies $(b,s)=(a,n)=1$.
(ii) The total major arc contribution with $(a,b)=1$ but $(b,s)>1$ or $(a,n)>1$ is $O_{A,k}(XY\mathcal{L}^{-A})$ for every fixed $A$.
(iii) This exceptional set in (ii) includes $bmn=ars$.
\end{lemma}
\begin{proof}
(i) For a prime $p\mid(b,s)$, we see $bmn-ars=a-b$ forces $p\mid a$. But this contradicts $(a,b)=1$, so we conclude $(b,s)=1$. Interchanging the endpoints gives $(a,n)=1$.

(ii) Fix $a,b$ and $p\mid b$. Set $k_0=mn$, $l_0=rs$. The relation $|bk_0-al_0|\ll Y$ has $O(1)$ choices of either coordinate when the other is fixed. By Cauchy--Schwarz, $\tau(pu)\le2\tau(u)$, and \eqref{vanu:divisor-moments},
\begin{align*}
  \sum_{\substack{k_0,l_0\le4X\\|bk_0-al_0|\ll Y,\ p\mid l_0}}\tau(k_0)\tau(l_0) &\ll\left(\sum_{k_0\le4X}\tau(k_0)^2\right)^{1/2} \left(\sum_{u\le4X/p}\tau(pu)^2\right)^{1/2} \ll Xp^{-1/2}\mathcal{L}^3.
\end{align*}
There are $O(Y^2)$ pairs $a,b$ and $k$ choices of $p\mid b$. Multiplying by the coefficient and kernel bounds gives
\begin{align*}
  \ll_k Y^2\mathcal{L}^{4C}\frac{\mathcal{L}^{3A_0}}Y\,Xp_{\min}^{-1/2}\mathcal{L}^3 \ll_{A,k}XY\mathcal{L}^{-A}.
\end{align*}
Reverse the endpoints for $(a,n)>1$.

(iii) If $bmn=ars$ with $(a,b)=1$, then $P^+(b)<z<P^-(r)$ forces $b\mid s$. Hence $(b,s)=b>1$.
\end{proof}
\paragraph{Absolute sums and the restricted form.}
For fixed $a,b\asymp Y$, \eqref{eq:014} and Cauchy--Schwarz give
\begin{align*}
  \sum_{\substack{aq-1\asymp X\\bq-1\asymp X}}\tau(aq-1)\tau(bq-1)&\ll\frac XY\mathcal{L}^3,\\
  \sum_{\substack{k_0,l_0\le4X\\|bk_0-al_0|\ll Y}}\tau(k_0)\tau(l_0)&\ll X\mathcal{L}^3.
\end{align*}
The second bound uses the same bounded coordinate fibres as above. Let $Q_Y^{\mathrm{raw,\,abs}}$ and $Q_Y^{\mathrm{maj,\,abs}}$ be the termwise absolute sums defining $Q_Y$ and $Q_Y^{\rm maj}$, respectively. Thus
\begin{align}
  Q_Y^{\mathrm{raw,\,abs}}&\ll_k Y^2\mathcal{L}^{4C}\frac XY\mathcal{L}^3 =XY\mathcal{L}^{4C+3},\label{eq:031}\\
  Q_Y^{\mathrm{maj,\,abs}}&\ll_k Y^2\mathcal{L}^{4C}\frac{\mathcal{L}^{3A_0}}Y\,X\mathcal{L}^3 =XY\mathcal{L}^{4C+3+3A_0}.\label{eq:032}
\end{align}
For $(a,b)>1$, the fixed-pair bound and \eqref{eq:abcoll} give total major arc contribution
\begin{align*}
  \ll_k Y^2e^{-\mathcal{L}^\alpha}\frac XY\mathcal{L}^{4C+3+3A_0} \ll_{E,k}XY\mathcal{L}^{-E}.
\end{align*}
Combining Lemmas~\ref{lem:shared} and~\ref{lem:exceptional-endpoints} with $Q_Y^{\rm min}=Q_Y-Q_Y^{\rm maj}$ gives
\begin{equation}  \label{eq:minor-restricted}
  \begin{aligned}
  Q_Y^{\rm min}={}&V^{-2}\sum_{\substack{a,b\in\mathcal{A},\ m,n,r,s\\(a,b)=1,\ (b,s)=(a,n)=1}}\eta_Y(a)\eta_Y(b)\alpha_m\overline{\alpha_r}\beta_n\overline{\beta_s}\\
  &\hspace{4em}\times H^{\rm min}(bmn-ars;a,b)+O_{E,k}(XY\mathcal{L}^{-E}).
  \end{aligned}
\end{equation}
All coefficient supports remain as in Proposition~\ref{prop:min}.

\subsection{Prime-list operators}  \label{sec:001}  \label{sec:input}
We use the prime-list construction of \cite[Section~4.2]{LWW}. Set
\begin{equation}  \label{eq:mompar}
  J = \lfloor \mathcal{L}^{\jmath} \rfloor.
\end{equation}
Choose ordered auxiliary primes $p_{i,1},\ldots,p_{i,J} \in \mathcal{P}_i$ and put
\begin{align*}
  d = \prod_{i=1}^k \prod_{l=1}^J p_{i,l}, \qquad D \le d < 2D, \qquad \mu^2(d) = 1, \qquad (d,abns) = 1.
\end{align*}
Write $\mathsf{v}=\binom{\mathsf{x}}{\mathsf{y}}$ and $\mathsf{v}'=\binom{\mathsf{x}'}{\mathsf{y}'}$, so that $\det(\mathsf{v},\mathsf{v}') = \mathsf{x}\mathsf{y}' - \mathsf{y}\mathsf{x}'$. For each $d\sim D$, we substitute in
\begin{equation}  \label{eq:lifted-endpoints}
  \mathsf{v} = \binom{bdm}{s}, \qquad \mathsf{v}' = \binom{adr}{n}, \qquad \frac{\det(\mathsf{v},\mathsf{v}')}{d} = bmn - ars = t.
\end{equation}
 On the support in \eqref{eq:minor-restricted}, Lemma~\ref{lem:exceptional-endpoints} gives $0 < |t| \ll Y = o(z)$. Since $(bd,s) = (ad,n) = 1$,
\begin{align*}
  p \mid \gcd(bdm,s) \quad\Longrightarrow\quad p \mid m,\quad p > z,\quad p \mid t,
\end{align*}
which is impossible; the same argument applies to $\mathsf{v}'$. Thus both endpoints are primitive.
\paragraph{States and coefficients.}
Write $|\mathsf{v}| = \max(|\mathsf{x}|,|\mathsf{y}|)$ and set
\begin{equation}  \label{eq:physbox}
  B = [DYM,16DYM] \times [N,2N].
\end{equation}
Set $B^* = \{\mathsf{v} \in B \cap \mathbb Z^2:(\mathsf{x},\mathsf{y}) = 1\}$. A state $(\mathsf{v},c)$ consists of $\mathsf{v} \in B^*$ and a divisor $c\mid\mathsf{x}$ formed from $J+1$ distinct primes in each $\mathcal{P}_i$. The vectors in \eqref{eq:lifted-endpoints} lie in $B^*$. Define
\begin{equation}  \label{eq:endpoint-function}
  f(\mathsf{v},c) = \overline{\alpha_{\mathsf{x}/c}}\beta_{\mathsf{y}}, \qquad |f(\mathsf{v},c)| \le \mathcal{L}^{2C}.
\end{equation}
The coefficients are zero outside their original supports. If $f(\mathsf{v},c) \ne 0$, then $\mathsf{x}/c$ is $z$-rough, so $c$ contains every active prime factor of $\mathsf{x}$, exactly once. Hence $f(\mathsf{v},c) = f(\mathsf{v})$, with $f(\mathsf{v})=0$ when no state lies above $\mathsf{v}$. In \eqref{eq:lifted-endpoints},
\begin{align*}
  f(\mathsf{v}) = \overline{\alpha_m}\beta_s, \qquad f(\mathsf{v}') = \overline{\alpha_r}\beta_n, \qquad \overline{f(\mathsf{v})}f(\mathsf{v}') = \alpha_m\overline{\alpha_r}\beta_n\overline{\beta_s}.
\end{align*}
For the fixed active family, every state has weight $\kappa_J=((J+1)!)^kV^{-(J+1)}$. Define the inner product
\begin{equation}  \label{eq:space}
  \langle h_1,h_2\rangle = \kappa_J\sum_{(\mathsf{v},c)\in\mathscr S} \overline{h_1(\mathsf{v},c)}h_2(\mathsf{v},c)
\end{equation}
on the Hilbert space $\mathcal{H} = \ell^2(\mathscr S,\kappa_J)$.
We suppress the dependence on $D$ within one range.
\paragraph{The transition and its pairing.}
Average over the $(J+1)^k$ divisors $b\mid c$ in $\mathcal A$, set $d=c/b$, and replace $b$ by $a\in\mathcal A$ with $(a,c)=1$. Define
\begin{equation}  \label{eq:transition}
  \begin{aligned}
  (\mathsf{T}h)(\mathsf{v},c) ={}& \frac{1}{(J+1)^k}\sum_{\substack{b\mid c,\ b\in\mathcal A\\D\le d<2D}}\frac{D}{d}V^{-1}
    \sum_{\substack{a\in\mathcal A\\(a,c)=1}}\ \sum_{\mathsf{v}':(\mathsf{v}',ad)\in\mathscr S}\\
  &\hspace{2em}\times\eta_Y(a)\eta_Y(b)H^{\rm min}\!\left(\frac{\det(\mathsf{v},\mathsf{v}')}{d};a,b\right)h(\mathsf{v}',ad).
  \end{aligned}
\end{equation}
An allowed transition is an edge; a sequence of edges is a path, and its initial vector is its base point. The target sum is unweighted. The inner product contributes one factor $\kappa_J$, giving
\begin{equation}  \label{eq:padmass}
  \frac1D\frac{\kappa_J}{(J+1)^k}V^{-1}\frac D{d} = V^{-2}\frac{(J!)^k}{d V^J}.
\end{equation}
Expanding the pairings and summing the dyadic ranges gives the exact identity
\begin{equation}  \label{eq:allpads}
  \begin{aligned}
  \sum_{D}\frac{\langle f_{D},\mathsf{T}_{D}f_{D}\rangle}{D} &= V^{-2}\sum_{\substack{a,b\in\mathcal{A},\ m,n,r,s\\(a,b)=1,\ (b,s)=(a,n)=1}}\eta_Y(a)\eta_Y(b)\alpha_m\overline{\alpha_r}\beta_n\overline{\beta_s}\\
  &\qquad\quad\times H^{\rm min}(bmn-ars;a,b)\sum_{\substack{p_{i,l}\in\mathcal{P}_i\ (1\le l\le J)\\\mu^2(d)=1,\ (d,abns)=1}}\frac1{dV^J}.
  \end{aligned}
\end{equation}
For each outer summand, each admissible $d$ determines one edge and exactly one dyadic range, and occurs $(J!)^k$ times in the auxiliary sum. The condition $(d,ns)=1$ is required because $n,s$ need not be rough.
\subsection{Averaging the auxiliary primes}  \label{sec:pads}
\begin{lemma}[Transition formulation]\label{lem:transition-reconstruction}
Under the preceding hypotheses, for every fixed $E>0$,
\begin{equation}\label{eq:transition-reconstruction}
  Q_Y^{\rm min} = \sum_{D}\frac{\langle f_{D},\mathsf{T}_{D}f_{D}\rangle}{D} + O_{E,k}(XY\mathcal{L}^{-E}).
\end{equation}
\end{lemma}
\begin{proof}
The auxiliary mass is
\begin{align*}
  \sum_{\substack{p_{i,l}\in\mathcal{P}_i\\1\le l\le J}}\frac1{dV^J} = V^{-J}\prod_i\left(\sum_{p\in\mathcal{P}_i}\frac1p\right)^J = 1.
\end{align*}
Since $ns \le x^{O(1)}$ has $O(\mathcal{L})$ distinct prime divisors, a union bound gives
\begin{equation}  \label{eq:030}
  \begin{aligned}
  &0 \le 1 - \sum_{\substack{p_{i,l}\in\mathcal{P}_i\ (1\le l\le J)\\\mu^2(d)=1,\ (d,abns)=1}}\frac1{dV^J}\\
  &\qquad\ll_k \bigl((J+1)^2+(J+1)\mathcal{L}\bigr)\max_{i,p}(pV_i)^{-1} \ll_k \bigl((J+1)^2+(J+1)\mathcal{L}\bigr)e^{-\mathcal{L}^\alpha}.
  \end{aligned}
\end{equation}
For integer $t$ and $a,b$ on the support of $\eta_Y(a)\eta_Y(b)$,
\begin{align*}
  |H^{\rm min}(t;a,b)| \le \mathbf1_{t=a-b}+|H^{\rm maj}(t;a,b)|.
\end{align*}
For $(a,b)=1$, the raw term is the common-cofactor condition. Subtracting \eqref{eq:allpads} from \eqref{eq:minor-restricted} and using \eqref{eq:031}--\eqref{eq:032},
\begin{align*}
  &\left|Q_Y^{\rm min}-\sum_{D}\frac{\langle f_{D},\mathsf{T}_{D}f_{D}\rangle}{D}\right| \ll_{E,k}XY\mathcal{L}^{-E}+XY\mathcal{L}^{4C+3+3A_0}\bigl((J+1)^2+(J+1)\mathcal{L}\bigr)e^{-\mathcal{L}^\alpha}
\end{align*}
which is $O_{E,k}(XY\mathcal{L}^{-E})$.
\end{proof}
\subsection{Good states}  \label{sec:good-states}
Fix $D$. For primitive $\mathsf{v}=\binom{\mathsf{x}}{\mathsf{y}}$ with $\mathsf{x}>0$, choose $\mathsf{x}_0,\mathsf{y}_0$ with $\mathsf{x}\mathsf{y}_0-\mathsf{y}\mathsf{x}_0=-1$ and set $r_{\mathsf{v}}=\mathsf{x}_0/\mathsf{x}+\mathbb Z\in\mathbb T$. This is independent of the completion. We use the tests in \cite[(4.14)--(4.16)]{LWW}; in each, $d$ ranges over $c/b$ with $b\mid c$ and $b\in\mathcal A$.
\begin{enumerate}
\item For every such $d\in[D,2D)$,
\begin{equation}  \label{eq:goodone}
  \|l dr_{\mathsf{v}}\| > Y^{-0.7} \qquad(l\in\mathbb Z,\ 1\le l\le Y^{0.2}).
\end{equation}
\item Fix nonempty $I\subseteq\{1,\ldots,k\}$, set $i_*=\max I$, and let $Z_{\mathrm{prod}}=\prod_{i\in I\setminus\{i_*\}}q_i$, with $q_i\in\mathcal{P}_i$. Write $d',Z_{\mathrm{prod}}'$ for a second choice. Give $T_f=\prod_{i\notin I}p_i$, with $p_i\in\mathcal{P}_i$, the weight $\prod_{i\notin I}(p_iV_i)^{-1}$. Except for choices of $T_f$ of total weight at most $e^{-\lambda_{i_*}/4}$, require
\begin{equation}  \label{eq:goodtwo}
  \|(dZ_{\mathrm{prod}}-d'Z_{\mathrm{prod}}')T_fr_{\mathsf{v}}\| \ge 100e^{-\lambda_{i_*}} \qquad(dZ_{\mathrm{prod}}\ne d'Z_{\mathrm{prod}}').
\end{equation}
\end{enumerate}
For fixed $c$ and uniform $r\in\mathbb T$, the first failure set has measure at most $2(J+1)^kY^{-1/2}$. For the second,
\begin{align*}
  \#\{dZ_{\mathrm{prod}}\} &\le (J+1)^k\exp\!\left(\frac{2\lambda_{i_*}}{r_{\rm sep}-1}\right),\\
  \operatorname{meas}\{\text{the second test fails for }I\} &\le 200e^{\lambda_{i_*}/4-\lambda_{i_*}}\bigl(\#\{dZ_{\mathrm{prod}}\}\bigr)^2 \le 200(J+1)^{2k}e^{-\gamma_{\mathrm{sep}}\lambda_{i_*}},\\
  \gamma_{\mathrm{sep}} &= \frac34-\frac4{r_{\rm sep}-1}>0.
\end{align*}
Here the second line is Markov's inequality after averaging over $T_f$. Since $\lambda_{i_*}\ge\mathcal{L}^\alpha$ and $\log Y\ge\lambda_1-\log4$,
\begin{equation}  \label{rev:good-state-uniform}
  \operatorname{meas}\{r:\text{the state is not good}\} \ll_k e^{-c_{\mathrm{sep}}\mathcal{L}^\alpha}, \qquad c_{\mathrm{sep}}=\min\{1/8,\gamma_{\mathrm{sep}}/2\}>0.
\end{equation}
Let $\mathsf{G}$ multiply by the good-state indicator. Define
\begin{equation}  \label{eq:aone}
  \begin{aligned}
  \mathsf{A}&=\mathsf{G}\mathsf{T}\mathsf{G},\\
  \langle f,\mathsf{A}f\rangle&=\langle\mathsf{G}f,\mathsf{T}\mathsf{G}f\rangle.
  \end{aligned}
\end{equation}
Consequently,
\begin{equation}  \label{eq:bad-decomposition}
  \langle f,\mathsf{T}f\rangle-\langle f,\mathsf{A}f\rangle = \langle(1-\mathsf{G})f,\mathsf{T}f\rangle+\langle\mathsf{G}f,\mathsf{T}(1-\mathsf{G})f\rangle.
\end{equation}

\subsection{Reduction to an averaged moment}  \label{sec:minor-completion}  \label{sec:ends}
\begin{lemma}[Operator formulation]\label{lem:operator-reconstruction}
Under the preceding hypotheses, for every fixed $E>0$,
\begin{equation}  \label{eq:operator-reconstruction}
  Q_Y^{\rm min}=\sum_{D}\frac1{D}\langle f_{D},\mathsf{A}_{D}f_{D}\rangle_{\mathcal{H}_{D}}+O_{E,k}(XY\mathcal{L}^{-E}).
\end{equation}
\end{lemma}
The proof is given in Section~\ref{sec:minor-equidistribution}.
\paragraph{Endpoint norms.}
For a fixed admissible product $c$, there are $O(DYMN/c)$ possible vectors. Hence
\begin{align*}
  \kappa_J\sum_{c\text{ admissible}}\frac1{c}\le1,\\
  \|f\|^2\ll DYMN\mathcal{L}^{4C}\kappa_J\sum_{c\text{ admissible}}\frac1{c}.
\end{align*}
Thus
\begin{equation}  \label{eq:endnorm}
  \|f\|\ll(DYMN)^{1/2}\mathcal{L}^{2C},\qquad \sup|f|\le\mathcal{L}^{2C}.
\end{equation}
Set
\begin{equation}  \label{eq:moment-order}
  \begin{aligned}
  R&=\lceil\mathcal{L}^{\varrho}/2\rceil,\qquad u_{\mathsf{v}}(\mathsf{v}',c)=\mathbf1_{\mathsf{v}'=\mathsf{v}},\\
  \mathcal{T}_R&=\sum_{\mathsf{v}\in B^*}\langle u_{\mathsf{v}},\mathsf{A}^{2R}u_{\mathsf{v}}\rangle.
  \end{aligned}
\end{equation}
\begin{lemma}[Interpolation by endpoint counting]\label{lem:block-interpolation}
For the endpoint vector above,
\begin{equation}  \label{eq:028}
  |\langle f,\mathsf{A}f\rangle| \ll\|f\|^{2-1/R}\bigl(DY(\sup|f|)^2\mathcal{T}_R\bigr)^{1/(2R)}.
\end{equation}
\end{lemma}
\begin{proof}
Set $S=DY$ and write $\sigma(\mathsf{v})=\mathsf{y}/\mathsf{x}$. On each edge, the cutoff in \eqref{eq:arcs-kernels} gives $|\det(\mathsf{v},\mathsf{v}')|/d\ll Y$ in \eqref{eq:transition}. Since $d<2D$ and $\mathsf{x},\mathsf{x}'\asymp SM$ by \eqref{eq:physbox},
\begin{align*}
  |\sigma(\mathsf{v}')-\sigma(\mathsf{v})|=\frac{|\det(\mathsf{v},\mathsf{v}')|}{\mathsf{x}\mathsf{x}'}\ll\frac1{SM^2}.
\end{align*}
Summing along a path of length $2R$ gives $|\sigma(\mathsf{v}')-\sigma(\mathsf{v})|\ll R/(SM^2)$ at its endpoints, hence $|\det(\mathsf{v},\mathsf{v}')|\ll RS$. For fixed primitive $\mathsf{v}$ and fixed determinant, the possible $\mathsf{v}'$ differ by integer multiples of $\mathsf{v}$. Their first coordinates lie in $[SM,16SM]$ and are spaced by at least $\mathsf{x}\ge SM$. Thus each of the $O(RS)$ integer determinants gives $O(1)$ endpoints.

By self-adjointness, $\langle\mathsf{A}^Ru_{\mathsf{v}},\mathsf{A}^Ru_{\mathsf{v}'}\rangle=\langle u_{\mathsf{v}},\mathsf{A}^{2R}u_{\mathsf{v}'}\rangle$ vanishes unless $|\det(\mathsf{v},\mathsf{v}')|\ll RS$. Write $x_{\mathsf{v}}=\|\mathsf{A}^Ru_{\mathsf{v}}\|$. Expanding $f=\sum_{\mathsf{v}} f(\mathsf{v})u_{\mathsf{v}}$ and applying Cauchy--Schwarz gives
\begin{align*}
  \|\mathsf{A}^Rf\|^2
  &\le (\sup|f|)^2\sum_{|\det(\mathsf{v},\mathsf{v}')|\ll RS}x_{\mathsf{v}}x_{\mathsf{v}'}\\
  &\le \tfrac12(\sup|f|)^2\sum_{|\det(\mathsf{v},\mathsf{v}')|\ll RS}(x_{\mathsf{v}}^2+x_{\mathsf{v}'}^2)
   \ll RS(\sup|f|)^2\mathcal{T}_R,
\end{align*}
since $\mathcal{T}_R=\sum_{\mathsf{v}} x_{\mathsf{v}}^2$ by \eqref{eq:moment-order}. H\"older's inequality for the spectral measure of $\mathsf{A}$, as in \cite[(4.64)--(4.67)]{LWW}, gives
\begin{align*}
  |\langle f,\mathsf{A}f\rangle|
  \le \|f\|^{2-1/R}\|\mathsf{A}^Rf\|^{1/R}.
\end{align*}
The claim follows since $R^{1/(2R)}\ll1$.
\end{proof}
Equations \eqref{eq:endnorm} and \eqref{eq:028} give
\begin{align*}
  \frac{|\langle f,\mathsf{A}f\rangle|}{DYMN}&\ll\mathcal{L}^{4C}(DY)^{1/(2R)}\left(\frac{\mathcal{T}_R}{DYMN}\right)^{1/(2R)},\\
  (DY)^{1/(2R)}&=\exp\!\left(O_k(\mathcal{L}^{\beta+\jmath-\varrho})\right)=O_k(1),
\end{align*}
since $\beta+\jmath<\varrho$.
\begin{proposition}[Averaged moment bound]\label{prop:par}
Fix $E_0,\delta > 0$ and assume \eqref{eq:comparison-scales}. Under the stated prime-interval and parameter conditions, choosing $A_0$ sufficiently large and then a sufficiently large fixed $k$ gives
\begin{equation}  \label{eq:moment}
  \frac{\mathcal{T}_R}{DYMN} \ll \mathcal{L}^{-2RE_0}.
\end{equation}
The constants are uniform in the individual interval scales and the lengths in \eqref{eq:comparison-scales}.
\end{proposition}
The proof is given in Section~\ref{sec:minor-moments}. Hence using \eqref{eq:endnorm}
\begin{equation}  \label{eq:pair}
  |\langle f,\mathsf{A}f\rangle|\ll_k DYMN\mathcal{L}^{-E_0+4C}.
\end{equation}

Choose $E_0>A+4C+2$ and then $A_0,k$ as in Proposition~\ref{prop:par}. Thus \eqref{eq:operator-reconstruction} and \eqref{eq:pair}, with $E=A+2$, give
\begin{align*}
  |Q_Y^{\rm min}| &\le\sum_{D}\frac{|\langle f_{D},\mathsf{A}_{D}f_{D}\rangle|}{D}+O_{A,k}(XY\mathcal{L}^{-A-2}) \ll_{A,k}XY\mathcal{L}^{-E_0+4C+1}+XY\mathcal{L}^{-A-2} \ll_A XY\mathcal{L}^{-A}.
\end{align*}

This completes the proof of the minor arc Proposition~\ref{prop:min}, and hence the final Theorem~\ref{thm:chowla}. Hence it remains to prove Lemma \ref{lem:operator-reconstruction} and Proposition \ref{prop:par}.

\section{Equidistribution and operator formulation}  \label{sec:minor-equidistribution}

\subsection{Congruence equidistribution}  \label{sec:roots}

Set $A = DYM$ and $B_0 = N$. For positive coprime $\mathsf{x},\mathsf{y}$ with $\mathsf{x} \ge 2$, choose $0 \le \mathsf{x}_0 < \mathsf{x}$ with $\mathsf{y}\mathsf{x}_0 \equiv1 \pmod \mathsf{x}$, and set $\mathsf{y}_0 = (\mathsf{y}\mathsf{x}_0 - 1)/\mathsf{x}$, $r = \mathsf{x}_0/\mathsf{x}$, and $\tau(\mathsf{w}) = w_1 + r w_2$ for $\mathsf{w}=\binom{w_1}{w_2}$. Then
\begin{equation}  \label{eq:gprime}
  G = \begin{pmatrix}\mathsf{x}&\mathsf{x}r\\\mathsf{y}&\mathsf{y}r - 1/\mathsf{x}\end{pmatrix}, \qquad \det G = -1, \qquad G\mathsf{w} = \binom{\mathsf{x}\tau(\mathsf{w})}{\mathsf{y}\tau(\mathsf{w}) - w_2/\mathsf{x}}.
\end{equation}
Here $r = \mathsf{x}_0/\mathsf{x}$ is the normalized modular inverse of $\mathsf{y}$ modulo $\mathsf{x}$. Since $G$ is an integer matrix of determinant $-1$, we have
\begin{align*}
  G\mathcal{C} = \mathcal{C}, \qquad G\binom{1}{0} = \binom{\mathsf{x}}{\mathsf{y}}, \qquad \det(G\mathsf{w},G\mathsf{w}') = -\det(\mathsf{w},\mathsf{w}').
\end{align*}
Thus $H^{\rm min}(-t;a,b) = H^{\rm min}(t;b,a)$ gives the relative kernel
\begin{align*}
  \frac{D}{d}\eta_Y(a)\eta_Y(b)H^{\rm min}\!\left(\frac{\det(\mathsf{w},\mathsf{w}')}{d};b,a\right).
\end{align*}

For a prime $p$, let $\ell_p(\mathsf{v_0})$ be the reduction of the integer line
\begin{align*}
  \ell=G^{-1}\mathbb Z\binom{0}{1}=\mathbb Z\binom{-\mathsf{x}_0}{\mathsf{x}}.
\end{align*}
Write $\bar{\mathsf{w}}$ for the reduction of $\mathsf{w}$ modulo the prime under consideration. Since $\det G=-1$, the first row of $G$ is nonzero modulo $p$. For $\mathsf{v}=G\mathsf{w}=\binom{\mathsf{x}_{\mathsf{v}}}{\mathsf{y}_{\mathsf{v}}}$, we have
\begin{align*}
  p\mid \mathsf{x}_{\mathsf{v}} \quad\Longleftrightarrow\quad \mathsf{x}w_1+\mathsf{x}_0w_2\equiv0\pmod p \quad\Longleftrightarrow\quad \bar{\mathsf{w}}\in\ell_p(\mathsf{v_0}).
\end{align*}
The $p+1$ lines through the origin in $\mathbb F_p^2$ form $\mathbb P^1(\mathbb F_p)$. For primitive integer vectors, both reductions are nonzero, so
\begin{align*}
  \bar{\mathsf{w}},\bar{\mathsf{w}}'\text{ lie on a common line}
  \quad\Longleftrightarrow\quad p\mid\det(\mathsf{w},\mathsf{w}').
\end{align*}

For a positive squarefree integer $S$, put
\begin{align*}
  \Gamma=\mathrm{SL}_2(\mathbb Z), \qquad \widetilde G=\begin{pmatrix}\mathsf{y}&\mathsf{y}_0\\\mathsf{x}&\mathsf{x}_0\end{pmatrix}\in\Gamma, \qquad \Gamma_0(S)=\{\gamma\in\Gamma:\gamma_{21}\equiv0\pmod S\}.
\end{align*}
For two normalized completions $G,G'$, equality of their prescribed lines is equivalent to proportionality of their first rows modulo each $p\mid S$. Explicitly,
\begin{align*}
  \ell_p(\mathsf{v_0})=\ell_p(G')\ (p\mid S) &\quad\Longleftrightarrow\quad \mathsf{x}\mathsf{x}_0'-\mathsf{x}_0\mathsf{x}'\equiv0\pmod S\\
  &\quad\Longleftrightarrow\quad\widetilde G(\widetilde G')^{-1}\in\Gamma_0(S)\quad\Longleftrightarrow\quad\Gamma_0(S)\widetilde G=\Gamma_0(S)\widetilde G'.
\end{align*}
Hence the line choices correspond to the cosets in $\Gamma_0(S)\backslash\Gamma$, and
\begin{align*}
  [\Gamma:\Gamma_0(S)]=\prod_{p\mid S}(p+1).
\end{align*}

\begin{corollary}[Equidistribution in congruence cosets]\label{cor:root-slopes}
Let $A,B_0\ge2$, let $S$ be a positive squarefree integer, and prescribe a line $\ell_p\in\mathbb P^1(\mathbb F_p)$ for each $p\mid S$. For intervals $I_1,I_2\subseteq[1,16]$ and $I_3\subseteq[0,1]$, we have
\begin{equation}  \label{eq:app-root-slopes}
  \begin{aligned}
  &\frac1{AB_0}\#\left\{(\mathsf{x},\mathsf{y}):\begin{gathered} \mathsf{x},\mathsf{y}\ge1,\quad(\mathsf{x},\mathsf{y})=1,\quad \mathsf{x}/A\in I_1,\quad \mathsf{y}/B_0\in I_2,\\
  \mathsf{x}_0/\mathsf{x}\in I_3,\quad\ell_p(\mathsf{v_0})=\ell_p \ \forall \, p\mid S
  \end{gathered}\right\}\\
  &\quad=\frac{|I_1||I_2||I_3|}{\zeta(2)}\prod_{p\mid S}\frac1{p+1}+O\!\left(S^5\min(A,B_0)^{-1/1000}\right),
  \end{aligned}
\end{equation}
with an absolute implied constant.
\end{corollary}
\begin{proof}
This follows from the joint equidistribution estimate \cite[Lemma~A.1]{LWW} with $G\left(\begin{smallmatrix}1&1\\0&-1\end{smallmatrix}\right)$, the normalized determinant-one completion with ratio $1-\mathsf{x}_0/\mathsf{x}$. A prescribed line $\ell_p(\mathsf{v_0})$ occurs in $p(p-1)$ of the $p(p^2-1)$ determinant-$(-1)$ matrix classes modulo $p$, giving density $1/(p+1)$. The Chinese remainder theorem gives the product of these densities. Summing the error over at most $S^3$ classes gives $O(S^5\min(A,B_0)^{-1/100})$.
\end{proof}

The main term is the product of Lebesgue measure in $(\mathsf{x}/A,\mathsf{y}/B_0,\mathsf{x}_0/\mathsf{x})$, the primitive-vector density $1/\zeta(2)$, and uniform measure on $\Gamma_0(S)\backslash\Gamma$.

The vectors $u_{\mathsf{v}}$ are unnormalized. In the ordinary matrix representation of $\mathsf A$, their expansion gives separate sums over the initial and final products:
\begin{align*}
  \langle u_{\mathsf{v}},\mathsf{A}^{2R} u_{\mathsf{v}}\rangle = \kappa_J\sum_{c:(\mathsf{v},c)\in\mathscr S}\sum_{c':(\mathsf{v},c')\in\mathscr S}(\mathsf{A}^{2R})_{(\mathsf{v},c),(\mathsf{v},c')}.
\end{align*}
The factor $\kappa_J$ occurs once, at the initial state. The relative vectors trace a loop of length $2R$, with $\mathsf{w}_0=\mathsf{w}_{2R}=\binom{1}{0}$.

Extend $G$ to real parameters $(\mathsf{x},\mathsf{y},r)\in[A,16A]\times[B_0,2B_0]\times[0,1]$. For primitive $\mathsf{w}$, choose an integer vector $\widehat{\mathsf{w}}$ with $\det(\mathsf{w},\widehat{\mathsf{w}})=1$ and $|\widehat{\mathsf{w}}|\ll1+|\mathsf{w}|$; then the transported completion parameter is
\begin{align*}
  r_{G\mathsf{w}}\equiv\frac{\tau(\widehat{\mathsf{w}})}{\tau(\mathsf{w})}\pmod1.
\end{align*}
Define $\mathcal{M}_R(\mathsf{x},\mathsf{y},r)$ by the relative-coordinate expansion above, with primitive integer vectors and the same transition weights, adjoints, position restrictions and good-state tests. For a relative path $(\mathsf{w}_j,c_j)_{j\le 2R}$, set $J_p=\{j\le 2R:p\mid c_j\}$. For each $j\in J_p$, replace the divisibility condition at $p$ by $\mathbf1_{\bar{\mathsf{w}}_j\in\ell_p}$. Average over independent uniform lines $\ell_p\in\mathbb P^1(\mathbb F_p)$ at distinct primes. Every active occurrence of a given prime uses the same line.

\begin{proposition}[Averaging over primitive vectors]\label{prop:root-average}
Under \eqref{eq:comparison-scales}, for every fixed $E>0$,
\begin{align*}
  \frac{\mathcal{T}_R}{AB_0}=\frac1{\zeta(2)}\int_1^{16}\int_1^2\int_0^1\mathcal{M}_R(As,B_0t_0,r)\,dr\,dt_0\,ds+O_E(\mathcal{L}^{-2ER}).
\end{align*}
The same averaging formula holds for absolute one-edge sums, treating the raw and major arc kernels separately and allowing either endpoint to fail goodness. After bounding endpoint coefficients by their suprema, the normalized averaging error is $O_E(\mathcal{L}^{-E})$.
\end{proposition}
\begin{proof}
Along a supported path with vectors $\mathsf{v}_j=G\mathsf{w}_j=\binom{x_j}{y_j}$, the determinant cutoff gives
\begin{align*}
  \left|\frac{y_{j+1}}{x_{j+1}}-\frac{y_j}{x_j}\right|=\frac{|\det(\mathsf{v}_j,\mathsf{v}_{j+1})|}{x_j x_{j+1}}\ll\frac1{DYM^2},
\end{align*}
and hence
\begin{align*}
  |\mathsf{w}_j|\ll RDY,\qquad \frac1{16}\le\tau(\mathsf{w}_j)\le16.
\end{align*}
For a fixed path, let $S=\prod_{p:J_p\ne\varnothing}p$ be the product of its distinct prime values. Let $\mathfrak W_R$ be one plus the sum of absolute path weights over these bounded relative vectors, omitting the indicators that depend on the base point. Counting relative vectors and prime choices gives
\begin{align}
  \log S&\ll_k R(J+1)\mathcal{L}^{\beta},\label{eq:root-modulus-budget}\\
  \log\mathfrak W_R&\ll_k R\log(RDY)+R(J+1)\bigl(\mathcal{L}^{\beta}+\log(J+1)\bigr).\label{eq:root-path-mass}
\end{align}
By \eqref{par:all},
\begin{align*}
  \log S+\log\mathfrak W_R\ll_k\mathcal{L}^{\varrho+\beta+\jmath+o(1)}=o(\mathcal{L}),\qquad S,\mathfrak W_R=x^{o(1)}.
\end{align*}
Since $\min(A,B_0)\ge x^\delta$, the normalized box discrepancy in Corollary~\ref{cor:root-slopes} satisfies, uniformly over these configurations,
\begin{align*}
  S^5\min(A,B_0)^{-1/1000}\ll x^{-c_\delta+o(1)},\qquad c_\delta=\delta/1000.
\end{align*}
The total absolute mass $\mathfrak W_R=x^{o(1)}$ preserves this power saving in the averaging argument of \cite[Proposition~4.4]{LWW}.

Writing $\mathsf{x}=As$ and $\mathsf{y}=B_0t_0$, the boundaries of the position conditions are
\begin{equation}  \label{eq:faces}
  s\tau(\mathsf{w})\in\{1,16\},\qquad t_0\tau(\mathsf{w})-\frac{w_2}{AB_0s}\in\{1,2\}.
\end{equation}
Their derivatives in $s$ and $t_0$, respectively, are $\tau(\mathsf{w})\in[1/16,16]$, while
\begin{equation}  \label{eq:rderiv}
  \frac{\partial}{\partial r}\frac{\tau(\widehat{\mathsf{w}})}{\tau(\mathsf{w})}=\frac1{\tau(\mathsf{w})^2}.
\end{equation}
Also $|w_2|/(AB_0)\ll R/X$, and all remaining derivatives and numbers of boundary pieces are $x^{o(1)}$. Enumerate the products, prime choices and integer translates in the two good-state tests. Their indicators, including the weighted exceptional-choice condition, are constant between the resulting boundary values of $r$.

Take a mesh of side length $h=x^{-c_\delta/8}$ as in (4.30)--(4.32) of \cite[Proposition~4.4]{LWW}, or (3.31)--(3.34) of \cite[Lemma~3.4]{PD}. There are $O(h^{-3})$ cubes, and the boundary cubes have total volume $h x^{o(1)}$. Summing their box discrepancies also bounds their arithmetic mass. Thus the equidistribution and boundary errors, including all absolute path weights, satisfy
\begin{equation}  \label{eq:root-grid-errors}
  \frac{\mathcal{E}_{\rm count}+\mathcal{E}_{\partial}}{AB_0}\ll h^{-3}x^{-c_\delta+o(1)}+h x^{o(1)} =x^{-5c_\delta/8+o(1)}+x^{-c_\delta/8+o(1)}=O_E(\mathcal{L}^{-2ER}).
\end{equation}
The remaining kernel factors depend only on the fixed relative determinants and primes, and $\varrho<1$ gives the last bound. For one edge, omit the return condition on the vector, use $|\mathsf{w}|\ll DY$, and allow all primitive targets. Failure of goodness changes only the already controlled indicators, giving the stated one-edge versions.
\end{proof}

We can now prove the operator formulation for $Q_Y^{\rm min}$.
\begin{proof}[Proof of Lemma~\ref{lem:operator-reconstruction}]
Apply the one-edge part of Proposition~\ref{prop:root-average} to \eqref{eq:bad-decomposition}. With base point $\binom{1}{0}$, the relative target $\mathsf{w}$ satisfies
\begin{align*}
  1/16\le\tau(\mathsf{w})=w_1+r w_2\le16,\qquad
  \begin{cases}
  w_2=d(b-a),&\text{for the raw kernel},\\
  w_2/d\in\mathbb Z,\quad|w_2/d|\ll Y,&\text{for the major arc kernel}.
  \end{cases}
\end{align*}
There are $O(1)$ possible $w_1$ for each $w_2$. The raw target count is $O(1)$, and \eqref{eq:uniform-kernel-size} gives major arc target mass $O(Y)\mathcal{L}^{3A_0}/Y=O(\mathcal{L}^{3A_0})$. The normalized prime mass is at most
\begin{align*}
  V^{-(J+2)}\prod_i\left(\sum_{p\in\mathcal{P}_i}\frac1{p+1}\right)^{J+2}\le1.
\end{align*}
Using $\sup|f|^2\le\mathcal{L}^{4C}$ and \eqref{rev:good-state-uniform}, and reversing the edge for a bad target, gives
\begin{equation}  \label{eq:baderr}
  |\langle f,\mathsf{T}f\rangle-\langle f,\mathsf{A}f\rangle| \ll_{E,k}DYMN\left(\mathcal{L}^{4C+3A_0}e^{-c_{\mathrm{sep}}\mathcal{L}^\alpha}+\mathcal{L}^{-E}\right).
\end{equation}
This is the one-edge comparison in \cite[(4.69)--(4.76)]{LWW}; see also \cite[(4.56)--(4.62)]{PD}. Take precision $E+2$ in \eqref{eq:baderr}. Since there are $O_k(\mathcal{L})$ dyadic ranges and $MN=X$, \eqref{eq:transition-reconstruction} gives \eqref{eq:operator-reconstruction}.
\end{proof}

\section{Prime lifetimes and moments}  \label{sec:minor-moments}

Fix real parameters $(\mathsf{x},\mathsf{y},r)$ in the averaged moment. A prime is born at first appearance, remains live through temporary absences, and dies at final departure. A state records the active product $c$ and inactive live prime--line pairs, all with distinct prime values.

Choose $C_{\mathrm{pos}}$ large enough for the path bound in Proposition~\ref{prop:root-average}, and give counting measure to
\begin{equation}  \label{eq:relset}
  \mathcal{Z}_0 = \left\{\mathsf{w} \in \mathcal{C}: |\mathsf{w}| \le 2C_{\mathrm{pos}}RDY,\quad DYM/\mathsf{x} \le \tau(\mathsf{w}) \le 16DYM/\mathsf{x}\right\}.
\end{equation}
The position projection is multiplication by
\begin{equation}  \label{eq:posproj}
  \chi(\mathsf{w}) = \mathbf1_{G\mathsf{w} \in B}.
\end{equation}
For primes $p \in \mathcal{P}_i$, put
\begin{equation}  \label{eq:nudef}
  \nu_i(p) = \frac1{V_i(p + 1)}.
\end{equation}
Uniformly in the interval scales,
\begin{equation}  \label{eq:numass}
  v_i = \sum_{p \in \mathcal{P}_i}\nu_i(p) \le 1, \qquad \Delta = \max_{i,p}\nu_i(p) \ll_k e^{-\lambda_1}.
\end{equation}
The active product $c$ has weight $\prod_i\bigl((J+1)!\prod_{p\mid c,\,p\in\mathcal P_i}\nu_i(p)\bigr)$. The inactive tuple $\mathfrak m_i = ((p_l,\ell_l))_{l=1}^{k_i}$, with $p_l \in \mathcal{P}_i$, has measure $\prod_l\nu_i(p_l)/k_i!$ and counting measure on the lines $\ell_l\in\mathbb P^1(\mathbb F_{p_l})$. Functions are symmetric in inactive pairs. This is the factorial normalization of \cite[Sections~4.2--4.3]{PD}; inactive pairs persist with weight one. Let $\mathsf{P}_{\mathrm{dis}}$ enforce distinct live prime values and impose $\sum_i k_i \le M_0 := \lceil\mathcal{L}^2\rceil$.

Use the relative geometric transition with position and good-state projections, and unequal unshared primes across the edge. Alternate the kernel and its transposed conjugate as in $(\mathsf{A}\mathsf{A}^*)^R$. For shared product $d$, departing product $b$, and arriving product $a=\prod_i p_i$, require $D \le d < 2D$ and the equivalent shared-prime conditions
\begin{align*}
  d \mid \det(\mathsf{w},\mathsf{w}') \quad\Longleftrightarrow\quad \bar{\mathsf{w}},\bar{\mathsf{w}}'\text{ lie on a common line for each }p\mid d.
\end{align*}
For $S,I \subseteq \{1,\ldots,k\}$, the departing prime $q_i\mid b$ in $\mathcal P_i$ goes inactive when $i\in S$ and dies otherwise. In the inactive case, append the pair consisting of $q_i$ and the line containing $\bar{\mathsf{w}}$ modulo $q_i$ to $\mathfrak m_i$. When $i\in I$, choose an inactive pair $(p,\ell)$ from the source, impose $\bar{\mathsf{w}}'\in\ell$, and remove the pair from $\mathfrak m_i$ as $p$ becomes active again. Otherwise a prime is born with measure $\nu_i$. A prime cannot go inactive and become active again on the same edge. The transfer factor is
\begin{align*}
  \frac{1}{(J+1)^k V}\prod_{i \notin I}\frac1{p_i + 1}.
\end{align*}
Summing over $S,I$ and the prime choices, and compressing on both sides by $\mathsf{P}_{\mathrm{dis}}$ and the cutoff, defines $\mathsf{E}_j$ for $0 \le j < 2R$.

For the set $J_p$ of vertices where $p$ is active, the average over lines in \cite[Section~4.1]{PD}, without damping, is
\begin{equation}  \label{eq:pfactor}
  \frac1{p + 1}\sum_{\ell\in\mathbb P^1(\mathbb F_p)}\prod_{j\in J_p}\mathbf1_{\bar{\mathsf{w}}_j\in\ell}.
\end{equation}
This equals one for $J_p=\varnothing$. Otherwise it equals $1/(p+1)$ if the reductions at the active vertices lie on a common line, and zero if they do not.

Initial entries are born at the first vertex. With distinct births and empty inactive tuples at both endpoints, each prime has one lifetime: birth, active runs, and death, if any. Primes active at the final vertex remain live. For $p \in \mathcal{P}_i$ with $\kappa_p$ active runs, the weight is
\begin{align*}
  \nu_i(p)V_i^{-(\kappa_p - 1)} = \frac{V_i^{-\kappa_p}}{p + 1},
\end{align*}
with the reductions at the active vertices on a common line. This is one factor $V_i^{-1}$ per run times \eqref{eq:pfactor}. Temporary absences preserve the live prime and its line. Deaths remove the prime, and distinct births introduce unused values. Hence $\mathsf{P}_{\mathrm{dis}}$ preserves each lifetime and its weight.

Let $b_0$ indicate $\mathsf{w} = \binom{1}{0}$ and empty inactive tuples, with the active product subject only to $\mathsf{P}_{\mathrm{dis}}$. Multiplying the primewise identities gives
\begin{equation}  \label{eq:undamped-history}
  \mathcal{M}_R(\mathsf{x},\mathsf{y},r) = \langle b_0,\mathsf{E}_0\cdots\mathsf{E}_{2R - 1}b_0\rangle_{\mathrm{birth\ distinct}}, \qquad \|b_0\|^2 \le \prod_i v_i^{J + 1} \le 1.
\end{equation}
At most $2kR < M_0$ primes can be inactive, so the cutoff preserves the identity. The distinct-birth restriction excludes rebirths of dead primes and is imposed in the proof of Proposition~\ref{prop:par} below.

The factorial calculation in (4.17)--(4.18) of \cite[Section~4.3]{PD} reverses the transfers and exchanges $S,I$, introducing the scalar
\begin{align*}
  \frac{\prod_{i \in S}V_i}{\prod_{i \in I}V_i} \asymp_k 1.
\end{align*}
Also $(\mathsf{P}_{\mathrm{dis}}\mathsf{E}\mathsf{P}_{\mathrm{dis}})^* = \mathsf{P}_{\mathrm{dis}}\mathsf{E}^*\mathsf{P}_{\mathrm{dis}}$, so the reversed source has distinct prime values. Write $R(\mathsf{E})$ for the supremum of the absolute row sums with respect to the state measure and $C(\mathsf{E}) = R(\mathsf{E}^*)$. Schur's test gives
\begin{align*}
  \|\mathsf{E}\| \le \sqrt{R(\mathsf{E})C(\mathsf{E})}.
\end{align*}
For fixed primitive $\mathsf{w}$, choose $\widehat{\mathsf{w}}$ with $\det(\mathsf{w},\widehat{\mathsf{w}}) = 1$. A target of determinant $j_1$ has the form
\begin{equation}  \label{eq:target}
  \mathsf{w}' = j_1\widehat{\mathsf{w}} + j_0\mathsf{w}, \qquad \frac{\tau(\mathsf{w}')}{\tau(\mathsf{w})} = j_0 + j_1\frac{\tau(\widehat{\mathsf{w}})}{\tau(\mathsf{w})} \in [1/16,16].
\end{equation}
There are $O(1)$ choices of $j_0$ for each $j_1$. The raw term fixes $j_1 = d(b - a)$; the major arc term has $O(Y)$ choices of $j_1/d$, each with kernel size $O(\mathcal{L}^{3A_0}/Y)$. Sum birth weights using $v_i \le 1$ and allow $M_0^k$ choices of inactive pairs, as in (4.23)--(4.25) of \cite[Section~4.4]{PD}. The same count for the adjoint gives
\begin{equation}  \label{eq:crudebd}
  R(\mathsf{E}_j),C(\mathsf{E}_j),\|\mathsf{E}_j\| \ll_k M_0^k(1 + \mathcal{L}^{3A_0}) \le \mathcal{L}^{c_{\mathrm{crude}}},
\end{equation}
where $c_{\mathrm{crude}} = c_{\mathrm{crude}}(k,A_0)$ is fixed. Bounded multipliers on prime choices preserve this bound if they are invariant under permutations of the source inactive tuples. Prescribed birth values contribute their $\nu_i$ weights to the row bound.

Write $\mathrm{unlink}$ for $S = I = \varnothing$ (unlinked transitions), and $\mathrm{link}$ otherwise (linked transitions). Use the minor arc kernel for $\mathrm{unlink}$, and bound the raw and major arc kernels separately for $\mathrm{link}$.

Partition $w_2/\tau(\mathsf{w})$ into intervals of length $DY$, and then by the prescribed lines at shared primes. The determinant cutoff gives bounded block adjacency. Fix two interacting intervals, the shared product $d$, and the inactive tuples. The line conditions at the shared primes define a lattice of index $d$. By the first good-state test and the lattice and Fourier argument in (4.42)--(4.46) of \cite[Section~4.4]{LWW}, there are a positively oriented basis $\mathcal{B}_{d}$ and a coordinate box $\mathcal{R}$ with side lengths
\begin{equation}  \label{eq:latbox}
  A'B' \asymp Y, \qquad Y^{1/10} \le A',B' \le Y^{9/10},
\end{equation}
such that
\begin{equation}  \label{eq:latbox-e}
  \mathsf{w} = \mathcal{B}_{d}\mathsf{k}, \quad \mathsf{w}' = \mathcal{B}_{d}\mathsf{l}, \quad \det\mathcal{B}_{d} = d, \quad \frac{\det(\mathsf{w},\mathsf{w}')}{d} = \det(\mathsf{k},\mathsf{l}).
\end{equation}
For $a = \prod_i p_i$, define $m(a) = V^{-1}\prod_{p\mid a}(p+1)^{-1}$; the disjoint intervals make this factorization unique. On the support of $\eta_Y(a)$,
\begin{align*}
  m(a) \ll_k Y^{-1}, \qquad \eta_Y(a)\sqrt{Ym(a)} \ll_k 1.
\end{align*}
After changing to counting measure, these factors and the position, goodness, primitivity, and prime-factor restrictions are bounded endpoint multipliers. With $d$ and the inactive tuples fixed, $\mathsf{P}_{\mathrm{dis}}$ also acts separately at the endpoints. Before excluding equal unshared primes across the edge, the enlarged Fourier operator in (4.45)--(4.46) of \cite[Section~4.4]{LWW} has norm
\begin{align*}
  O\!\left((\mathcal{L}^{-A_0} + Y^{-1/10})\log Y\right).
\end{align*}
This estimate precedes the specialization of $\log Y$ in the cited argument. Each excluded equality fixes a birth value of weight at most $\Delta$, so \eqref{eq:target} bounds both row and column sums of the collision operator by $O_k((1 + \mathcal{L}^{3A_0})\Delta)$. Subtracting this operator and using bounded block adjacency and Cauchy--Schwarz over $d$ (at most $(J+1)^k$ choices per endpoint) gives
\begin{equation}  \label{eq:cleanbd}
  \|\mathsf{E}_j^{\mathrm{unlink,min}}\| \ll_k (\mathcal{L}^{-A_0} + Y^{-1/10})\log Y + (1 + \mathcal{L}^{3A_0})\Delta \ll_k \mathcal{L}^{\beta - A_0} + Y^{-1/10}\log Y.
\end{equation}

Fix $S,I$ with $I\ne\varnothing$. Set $i_*=\max I$ and fix the returning pair $(p,\ell_p)$ with $p\in\mathcal{P}_{i_*}$. There are at most $M_0$ choices. Fix the newly born primes, of product $T_f=\prod_{i\notin I}p_i$, and set $U=ad/p$. Since $p\nmid c$, the raw determinant satisfies
\begin{equation}  \label{eq:rawdet}
  j_1=\det(\mathsf{w},\mathsf{w}')=bd-ad=c-pU
  \equiv c\not\equiv0\pmod p.
\end{equation}
Modulo $p$, the targets in \eqref{eq:target} lie on the affine line $c\overline{\widehat{\mathsf{w}}}+\mathbb{F}_p\bar{\mathsf{w}}$, which misses the origin since $\det(\mathsf{w},\widehat{\mathsf{w}})=1$. Thus $\ell_p$ admits at most one residue $j_0\equiv j_*\pmod p$, independent of $U$, as in \cite[(4.47)--(4.48)]{LWW} and \cite[(4.35)--(4.36)]{PD}. With $\rho=\tau(\widehat{\mathsf{w}})/\tau(\mathsf{w})$, the position restriction in \eqref{eq:target} gives
\begin{align*}
  U\rho\in\frac{j_*+c\rho-[1/16,16]}{p}\pmod1.
\end{align*}
This interval has length less than $16/p\le16e^{-\lambda_{i_*}}$. For nonexceptional $T_f$, \eqref{eq:goodtwo} separates the phases $U\rho$ of distinct $U$ by at least $100e^{-\lambda_{i_*}}$, since $U=dZ_{\mathrm{prod}}T_f$ and $\rho\equiv r_{G\mathsf{w}}\pmod1$. Hence at most one $U$ is admissible. As $(a,c)=1$, it determines
\begin{align*}
  d=\gcd(c,U),\qquad b=c/d,\qquad Z_{\mathrm{prod}}=U/(dT_f).
\end{align*}
Thus at most one of the $(J+1)^k$ choices of $b$ contributes, and $Z_{\mathrm{prod}}$ determines the other returning pairs by distinctness. The exceptional-weight bound applies since $\nu_i(p)\le(pV_i)^{-1}$. For exceptional birth choices use the crude $M_0^k$ count. The row bound for $I \ne \varnothing$ is at most
\begin{equation}  \label{eq:epsk}
  \eps_k = C_k\left(M_0(J + 1)^{-k} + M_0^ke^{-\gamma\lambda_1}\right)
\end{equation}
for a fixed $\gamma > 0$; for $I = \varnothing$ it is $O_k(1)$. The adjoint exchanges $S,I$ and preserves distinctness, giving the column bounds. If $\jmath k > 2$, then $\eps_k \le 1$ for sufficiently large $x$, and Schur's test gives, for every fixed $A$,
\begin{align*}
  \|\mathsf{E}_j^{\mathrm{link,raw}}\| \ll_k \sqrt{\eps_k} \ll_k \mathcal{L}^{1 - \jmath k/2} + O_A(\mathcal{L}^{-A}).
\end{align*}

For a prime becoming active again from a pair $(p,\ell)$, the condition $\bar{\mathsf{w}}'\in\ell$ is a nontrivial linear congruence (mod $p$). Distinctness implies $p \nmid d$, so $\mathcal{B}_{d}$ is invertible (mod $p$). By (4.40)--(4.41) of \cite[Section~4.7]{PD}, there are $O(Y(p^{-1} + Y^{-1/10}))$ admissible targets in $\mathcal{R}$. The kernel bound and the choices of inactive pairs bound the row by
\begin{align*}
  O_k\left(M_0^k\frac{\mathcal{L}^{3A_0}}Y\cdot Y(p^{-1} + Y^{-1/10})\right).
\end{align*}
Apply \eqref{eq:crudebd} in reverse if only $S$ is nonempty, flip the edge. By Schur's test, for fixed $c_{\mathrm{maj}} = c_{\mathrm{maj}}(k,A_0)$,
\begin{equation}  \label{eq:dirty-maj-bound}
  \|\mathsf{E}_j^{\mathrm{link,maj}}\| \ll_k \mathcal{L}^{c_{\mathrm{maj}}}(p_{\min}^{-1} + Y^{-1/10})^{1/2} \ll_A \mathcal{L}^{-A},
\end{equation}
since $p_{\min} \ge e^{\mathcal{L}^\alpha}$ and $\log Y \gg \mathcal{L}^\alpha$.

Given $E_1 > 0$, take $A_0 > E_1 + \beta + 2$ and then a fixed $k$ with $\jmath k/2 > E_1 + 3$. The decomposition
\begin{align*}
  \mathsf{E}_j = \mathsf{E}_j^{\mathrm{unlink,min}} + \mathsf{E}_j^{\mathrm{link,raw}} - \mathsf{E}_j^{\mathrm{link,maj}}
\end{align*}
and the three bounds give
\begin{equation}  \label{eq:esmall}
  \|\mathsf{E}_j\| \le \mathcal{L}^{-E_1} \qquad(0 \le j < 2R)
\end{equation}
for sufficiently large $x$, uniformly in the individual interval scales.

\begin{proof}[Proof of Proposition~\ref{prop:par}]
Take $E_1 = E_0 + 5$ and choose $A_0,k$ as in \eqref{eq:esmall}, fixing $c_{\mathrm{crude}}$ in \eqref{eq:crudebd}, enlarged if necessary. Apply the collision-rank argument in (4.44)--(4.55) of \cite[Section~4.8]{PD}. Every edge includes its current-state projections. Let $\mathcal{B}$ be the set of initial and possible later birth slots, with initial weights given by $\nu_i^{\otimes(J+1)}$ on distinct primes in each interval. Then
\begin{equation}  \label{eq:birthct}
  |\mathcal{B}| \le k(J + 1) + 2kR = \mathcal{L}^{O_k(1)}.
\end{equation}
An incoming prime is born or becomes active again. Birth measures have atoms at most $\Delta \ll_k e^{-\lambda_1}$. For $0\le s\le k$, keeping the weights of $s$ prescribed birth values explicit in \eqref{eq:crudebd} gives
\begin{equation}  \label{eq:eatom}
  \sup\{\text{absolute row sums with }s\text{ fixed birth values}\} \ll_k M_0^k(1 + \mathcal{L}^{3A_0})\Delta^s.
\end{equation}
Initial entries have the same atomic bound, and current-state exclusions decrease the absolute sums. The birth flags and Fourier phases in (4.53) of \cite[Section~4.8]{PD} depend only on birth choices and are invariant under permutations of the source inactive tuples. Thus each modified edge satisfies \eqref{eq:crudebd}, and each unmodified edge satisfies \eqref{eq:esmall}.

At collision rank $l$, the inclusion--exclusion terms have total coefficient mass at most $|\mathcal{B}|^{2l}$. High ranks gain $\Delta^l$ in absolute row iteration, and low ranks modify at most $2l$ edges. Set
\begin{equation}  \label{eq:rankcut}
  l_0 = \left\lceil C_0\frac{2R\log\mathcal{L}}{\lambda_1}\right\rceil, \qquad C_0 > 2(c_{\mathrm{crude}} + E_0 + 7).
\end{equation}
Since $\lambda_1 \ge \mathcal{L}^\alpha$, uniformly in the interval scales,
\begin{align*}
  C_k|\mathcal{B}|^2\Delta \le e^{-\lambda_1/2}, \qquad l_0\log\mathcal{L} = o(R).
\end{align*}
Thus the high and low rank estimates (4.52) and (4.55) of \cite[Section~4.8]{PD} imply
\begin{align*}
  \left|\langle b_0,\mathsf{E}_0\cdots\mathsf{E}_{2R - 1}b_0\rangle_{\mathrm{birth\ distinct}}\right| \ll \mathcal{L}^{-2R(E_0 + 4)}.
\end{align*}
By \eqref{eq:undamped-history} and Proposition~\ref{prop:root-average}, with any fixed precision $E > E_0$,
\begin{align*}
  \frac{\mathcal{T}_R}{DYMN} \ll \mathcal{L}^{-2R(E_0 + 4)} + \mathcal{L}^{-2ER} \ll \mathcal{L}^{-2RE_0}.
\end{align*}
\end{proof}

\begin{remark}[The reversed orientation]\label{det:reflected-orientation}
For $\mu(mn-1)$, the cofactor equation is $bmn-ars=b-a$ and the physical phase is $t-b+a$. Set
\begin{align*}
  G_+ = G\begin{pmatrix}1&1\\0&-1\end{pmatrix} =\begin{pmatrix}\mathsf{x}&\mathsf{x}-\mathsf{x}_0\\\mathsf{y}&\mathsf{y}-\mathsf{y}_0\end{pmatrix},\qquad \det G_+=1,\qquad r_+=1-r,\\
  \tau_+(\mathsf{w})=w_1+r_+w_2,\qquad G_+\mathsf{w}=\binom{\mathsf{x}\tau_+(\mathsf{w})}{\mathsf{y}\tau_+(\mathsf{w})+w_2/\mathsf{x}}.
\end{align*}
The relative equation remains $\det(\mathsf{w},\mathsf{w}')/d=b-a$. The fixed column operation bijects the matrix residue classes and the lines modulo $p$, and $r\mapsto1-r$ preserves uniform measure and both good-state tests. The sign change in the second coordinate preserves the boundary and operator bounds; primitivity uses only the absolute determinant bound. Thus the same estimates apply with $mn+1,rs+1$ replaced by $mn-1,rs-1$.
\end{remark}


\end{document}